\documentclass[11pt,reqno]{amsart}
\usepackage{verbatim}
\usepackage{amssymb}
\usepackage{enumerate}
\usepackage[active]{srcltx}
\numberwithin{equation}{section}       

\newtheorem{case}{Case}

\usepackage{xcolor}

\newif\ifdeveloping

\newcommand{\ical}{{\mathcal I}}

\newcommand{\tcal}{{\mathcal T}}

\newcommand{\setm}{\setminus}

\def\<{\left\langle}
\def\>{\right\rangle}
\def\cf{\operatorname{cf}}

\def\OO{{\omega}}
\def\br#1;#2;{\bigl[ {#1} \bigr]^ {#2} }

\def\ooseq#1;#2;{\< {#1}_{#2}:{#2}<\oo\>}
\def\ooset#1;#2;{\{ {#1}_{#2}:{#2}<\oo\}}
\def\seq#1;#2;#3;{\< {#1}_{#2}:{#2}<#3\>}
\def\set#1;#2;#3;{\{ {#1}_{#2}:{#2}<#3\}}
\def\oseq#1;#2;{\< {#1}_{#2}:{#2}<\OO\>}
\def\oset#1;#2;{\{ {#1}_{#2}:{#2}<\OO\}}
\def\oosequ#1;#2;{\< {#1}^{#2}:{#2}<\oo\>}
\def\oosetu#1;#2;{\{ {#1}^{#2}:{#2}<\oo\}}
\def\sequ#1;#2;#3;{\< {#1}^{#2}:{#2}<#3\>}
\def\setu#1;#2;#3;{\{ {#1}^{#2}:{#2}<#3\}}
\def\osequ#1;#2;{\< {#1}^{#2}:{#2}<\OO\>}
\def\osetu#1;#2;{\{ {#1}^{#2}:{#2}<\OO\}}

\newcommand{\al}{\alpha}
\newcommand{\be}{\beta}

\newcommand{\CS}{\operatorname{CS}}
\newcommand{\lop}{\operatorname{l}}
\newcommand{\kop}{\operatorname{k}}

\def\to{\longrightarrow}

\newcommand{\newcases}{\setcounter{case}{0}}

 \ifdeveloping
\usepackage[notref,notcite]{showkeys}
\fi
\usepackage{enumerate}

\newcommand{\prtime}{{\count0=\time\divide\count0 by 60
\count1=-\count0\multiply\count1 by 60 \advance\count1 by \time
\the\count0:\the\count1} }

\def\myheads#1;#2;{
\pagestyle{myheadings} \markboth{{\sc\hfill
#1\hfill\protect\makebox[0cm][r]{\rm\today; \prtime}}}
{{\sc\protect\makebox[0cm][l]{\rm\today;\ \prtime}\hfill
#2\hfill}} \thispagestyle{myheadings} }

\newcommand{\htt}{\operatorname{ht}}

\newcommand{\intpart}[1]{\operatorname{\mathcal E}({#1})}

\newcommand{\de}{\delta}
\newcommand{\ka}{\kappa}
\newcommand{\la}{\lambda}

\newtheorem{theorem}{Theorem}[section]
\newtheorem{proposition}[theorem]{Proposition}
\newtheorem{lemma}[theorem]{Lemma}

\newtheorem{claim}[theorem]{Claim}

\theoremstyle{definition}
\newtheorem{definition}[theorem]{Definition}

\newcommand{\ga}{\gamma}
\newcommand{\om}{\omega}
\def\<{\left\langle}
\def\>{\right\rangle}

\title{ A consistency theorem for cardinal sequences of length  $< \omega_3$}

\author[J. C. Martinez]{Juan Carlos Mart\'{\i}nez }
\address{Facultat de Matem\`atiques  i Inform\`atica\\ Universitat de Barcelona \\ Gran
 Via 585 \\ 08007 Barcelona, Spain}
 \email{jcmartinez@ub.edu}

\author[L. Soukup]{
Lajos Soukup }
\address{HUN-REN Alfr{\'e}d R{\'e}nyi Institute of Mathematics}
\email{soukup@renyi.hu}

\thanks{1.The research  and preparation of this paper was
supported by  NKFI grant K129211, the Spanish MICIN grant PID2020-116773GB-100 and the Catalan Goverment grant 2021SGR 00348.\\ \hspace*{3.5mm} 2. Corresponding author. \\ \hspace*{3.5mm}  {\em E-mail address}: jcmartinez@ub.edu (J.C. Mart\'{\i}nez) }

\begin{document}

\begin{abstract}
We prove that if $\lambda$ is a fixed uncountable cardinal and $f = \langle \ka_{\al} : \al < \delta \rangle$ is a sequence of infinite cardinals where $\delta < \omega_3$ and $\ka_{\al}\in \{\om,\lambda\}$ for each $\al < \delta$ in such a way that $f^{-1}\{\om\}$ is $\om_2$-closed in $\delta$, then it is consistent that there is a scattered Boolean space whose cardinal sequence is $f$.

\end{abstract}

\maketitle

\section {\bf Introduction}

If $X$ is a locally compact scattered Hausdorff (in short: LCS) space and $\al$ is an ordinal, the ${\alpha}^{th}$- {\em Cantor-Bendixson level} of $X$ is defined  by $I_{\al}(X) =$ the set of all isolated points of the subspace $X\setminus \bigcup\{I_{\beta}(X): \beta < \alpha \}$. The {\em height } of $X$ is defined as
$\htt(X) =$ the least ordinal $\de$ such that $I_{\de}(X) = \emptyset$, and the {\em reduced height} of $X$ is $\htt(X)^- =$ the least ordinal $\de$ such that $I_{\de}(X)$  is finite. Clearly, $\htt^-(X) \leq \htt(X) \leq \htt^-(X) + 1$. We  then define the {\em cardinal sequence} of $X$  as the sequence of the cardinalities of the infinite Cantor-Bendixson levels of $X$, i.e.
 $$\CS(X) = \langle |I_{\al}(X)| : \al < \htt^-(X)\rangle.$$

 We refer the reader to the survey papers \cite{ba} and \cite{ma2} for an extensive list of results on cardinal sequences of LCS spaces as well as examples and basic facts.
 It is well-known that any result on cardinal sequences of LCS spaces can be transferred to the context of superatomic Boolean algebras, since the one-point compactification of an LCS space is a scattered Boolean space,  and   the  class of  clopen algebras of scattered Boolean spaces coincides with the class of superatomic Boolean algebras.

\smallskip
  
  If $\lambda$ is an infinite cardinal and $\alpha$ is an ordinal, we denote by $\langle \lambda \rangle_{\alpha}$  the cardinal sequence $\langle \lambda_{\be} : \be < \al \rangle$ where $\lambda_{\be} = \lambda$ for all $\be < \al$.

\smallskip
\subsection*{Fundamental theorems and the main result}
 In \cite[Theorem 5]{jw}, Juh\'asz and Weiss characterized, in ZFC, all sequences of infinite cardinals of length $\om_1$ which arise as cardinal sequences of some LCS space. They  proved that
 if $f = \langle \ka_{\al} : \al < \omega_1 \rangle$ is a sequence
of infinite cardinals, then $f$ is the cardinal sequence of an LCS space iff $\ka_{\be}\leq \ka^{\omega}_{\al}$ for every $\al < \be < \omega_1$. It is known, however,  that this result can not be extended to cardinal sequences of length $\omega_1 + 1$ (see \cite[Theorem 3.2]{bs}). Nevertheless, it was shown in  \cite[Theorem 9]{jw} that if $f = \langle \ka_{\al} : \al < \delta \rangle$ is a sequence of infinite cardinals with $\delta < \om_2$ satisfying that $\ka_{\be}\leq \ka^{\omega}_{\al}$ for all $\al < \be < \delta$, and $\ka_{\al} \leq \om_1$ whenever $\cf (\al) = \om_1$, then there is an LCS space $X$ such that $\CS(X) = f$.  A characterization
under GCH for cardinal sequences of length $< \omega_2$ was later obtained  by Juh\'asz, Soukup and Weiss in \cite{jsw}.

\vspace{1mm} On the other hand, it was shown by Baumgartner and Shelah in \cite{bs} that it is relatively consistent with ZFC that there is an LCS space $X$ such that $\CS(X) = \langle \om \rangle_{\om_2}$. This result was improved by Soukup in \cite{so}, where it was shown that if GCH holds and $\la \geq \om_2$ is a regular cardinal, then in some cardinal-preserving generic extension we have $2^{\om} = \la$ and every sequence $f = \langle \ka_{\al} : \al < \om_2 \rangle$ of infinite cardinals with $\ka_{\al} \leq \la$  is the cardinal sequence of some LCS space. 

\vspace{1mm} Assume that $\delta$ is an ordinal with $\om_2 < \delta < \om_3$ and $L \subset \delta$. We say that $L$ is $\om_2$-{\em closed} in $\delta$ if for every strictly increasing sequence $\langle \al_{\xi} : \xi < \om_2 \rangle\in {}^{\om_2} L$, we have $\sup\{\al_{\xi} : \xi < \om_2 \}\in L \cup \{\delta\}$. 

\vspace{1mm} As a consequence of \cite[Proposition 1.1]{ms}, we obtain the following fact.

\begin{proposition}
If  $\delta$ is an ordinal with $\om_2 < \delta < \om_3$ and $f = \langle \ka_{\xi} : \xi < \delta \rangle$ is the cardinal sequence of some LCS space, then  $f^{-1}\{\om\}$ is  $\om_2$-closed in $\delta$.
\end{proposition}

In particular, a cardinal sequence $\langle \ka_{\xi} : \xi \leq \om_2 \rangle$ where $\ka_{\xi} = \om$ for $\xi < \om_2$ and $\ka_{\om_2} > \om$ can not be the cardinal sequence of an LCS space.  So,  Soukup's theorem can not be extended to cardinal sequences of length $< \om_3$.

\vspace{1mm}
However, it was proved in \cite{ms} that for any ordinal ${\omega}_2<{\delta}<{\omega}_3$
it is consistent that many sequences of cardinals of length ${\delta}$
are cardinal sequences of LCS spaces:


\begin{theorem}[{\cite[Theorem~2.1]{ms}}]\label{tm:8main}
  If GCH holds and $\lambda \ge \omega_2$ is a regular cardinal, then in some
  cardinal-preserving generic extension we have $2^{\omega} = \lambda$, and for every ordinal 
  $\eta < \omega_3$
  and every sequence $f = \langle \kappa_{\alpha} : \alpha < \eta \rangle$ of infinite cardinals with 
  $\kappa_{\alpha} \le \lambda$ for all $\alpha < \eta$ and 
  $\kappa_{\alpha} = \omega$ whenever $\cf(\alpha) = \omega_2$, the sequence $f$ is the cardinal sequence of
  some LCS space.
  \end{theorem}

In this paper we will prove the following result.

\begin{theorem}\label{Thm.1.2} 
  Assume that $GCH + \square_{\omega_1}$ holds and that $\lambda \geq \om_2$ is a regular cardinal. Let $\delta < \om_3$ and let $f = \langle \ka_{\al} : \al < \delta \rangle$ be a sequence of infinite cardinals such that $\ka_{\al}\in \{\om,\lambda \}$ for each $\al < \delta$ and such that $f^{-1}\{\om\}$ is $\om_2$-closed in $\de$. Then, in some cardinal-preserving generic extension with  $2^{\om} = \la$, there exists an LCS space $X$ with $\CS(X) = f$.
\end{theorem}

Moreover, the space $X$ constructed in Theorem \ref{Thm.1.2} will satisfy  a stronger property: for every uncountable cardinal $\lambda' \le \lambda$ there is an open subspace $X'$ of $X$ such that $\CS(X') = \langle \lambda_{\al} : \al < \de \rangle$ where $\lambda_{\al} = \om$ if $\ka_{\al} = \om$ and $\lambda_{\al} = \lambda'$ if $\ka_{\al} = \lambda$.

Theorem  \ref{tm:8main} will be used in the proof of Theorem \ref{Thm.1.2}, but  
these two theorems involve incomparable restrictions
concerning cardinal sequences, and thus neither of them implies the other.  Also, the proof of Theorem 1.3 is much more involved than the proof of Theorem 1.2.

However, if $f = \langle \ka_{\al} : \al < \delta \rangle$ is a fixed cardinal sequence where $\om_2 < \delta < \om_3$ such that each $\ka_{\al} \in \{\omega,\omega_1,\omega_2\}$, we do not know whether it is consistent that  there is an LCS space whose cardinal sequence is $f$. And we do not know it, even if $\delta = \om_2 + 1$.

Next,   we recall some notions and results
which will be useful in the proof of Theorem~\ref{Thm.1.2}, 
beginning with the combinatorial notion of a $\Delta$-function on $\omega_2$ introduced by Shelah
in \cite{bs}.

\begin{definition}
Let $A$ be a set of ordinals of order type $\om_2$. Let $F : [A]^2 \rightarrow [A]^{\leq \om}$ be such that $F\{\al,\be\}\subset \mbox{min}\{\al,\be\}$ for all $\{\al,\be\}\in [A]^2$.

\vspace{1mm}\noindent (a) We say that a set $E$ of finite subsets of $A$ is {\em adequate} for $F$, if for all $\{a,b\}\in [E]^2$, all $\al\in a\setminus b$, $\be \in b\setminus a$ and $\tau\in a \cap b$ the following holds:

\begin{enumerate}[(1)]
\item if $\tau < \al,\be$ then $\tau \in F\{\al,\be\}$,

\vspace{1mm}
\item if $\tau < \be$ then $F\{\al,\tau\}\subset F\{\al,\be\}$,

\vspace{1mm}
\item if $\tau < \al$ then $F\{\be,\tau\}\subset F\{\al,\be\}$.

\end{enumerate}

\vspace{2mm}\noindent (b) We say that $F$ is a $\Delta$-{\em function}, if for every uncountable set $D$ of finite subsets of $A$ there is an uncountable subset $E$ of $D$ that is adequate for $F$.

\end{definition}

The following well-known result is due to Veli\v{c}kovi\'c  (see \cite[Chapter 7 and Lemma 7.4.9]{to2} for a proof).

\begin{theorem}\label{Thm.1.4}
    $\square_{\om_1}$ implies the existence of a $\Delta$-function on $\om_2$.
\end{theorem}

Clearly, a $\Delta$-function on $\om_2$ can be transferred to a $\Delta$-function on any set of ordinals of order type $\om_2$.

We will also use the following notion, which permits us to construct in a direct way LCS spaces from partial orders.

If $ T = \bigcup\{\{\al\}\times A_{\al} : \al < \delta \}$ where $\delta$ is a non-zero ordinal and each $A_{\al}$ is a non-empty set of ordinals,  then for every $s = \langle \al,\zeta\rangle\in T$ we write $\pi(s)= \al$.

\begin{definition}\label{Df.1.7}
 We say that $\tcal = \langle T,\preceq,i \rangle$ is an  {\em LCS poset} on $T$, if the following conditions hold:

\begin{enumerate}
\item $\langle T,\preceq \rangle$ is a partial order with $T = \bigcup \{T_{\al}: \al < \delta \}$ for some non-zero ordinal $\delta$
such that each $T_{\al} = \{\al\}\times A_{\al}$ where $A_{\al}$ is a non-empty set of ordinals.

\item  If $s \prec t$ then $\pi(s) < \pi(t)$.

\item If $\al < \be < \eta$ and $t\in T_{\be}$, then $\{s\in T_{\al} : s \prec t \}$ is infinite.

\item   $i : [T]^2 \rightarrow [T]^{<\om}$ such that for every $\{s,t\}\in [T]^2$ the following holds:

\begin{enumerate}[(a)]
\item If $v\in i\{s,t\}$ then $v\preceq s,t$.

\item If $u\preceq s,t$ then there is  $v\in i\{s,t\}$
such that $u\preceq v$.
\end{enumerate}

\end{enumerate}
\end{definition}

\vspace{1mm}
If $\tcal = \langle T,\preceq,i \rangle$ is an LCS poset with $T = \bigcup \{T_{\al}: \al < \delta \}$, we define its {\em associated LCS space} $X = X(\tcal)$ as follows. The underlying set of $X(\tcal)$ is $T$. If $x\in T$, we write  $C(x)= \{y\in T : y\preceq x\}$. Then, for every $x\in T$ we define a basic neighborhood of $x$ in X as a set of the form $C(x)\setminus (C(x_1)\cup \dots \cup C(x_n))$ where $n < \omega$ and $x_1,\dots,x_n \prec x$. It can be checked that $X$ is a locally compact, scattered, Hausdorff space of height $\delta$ such that $I_{\al}(X) = T_{\al}$ for every $\al < \delta$ (see \cite{ba} for a proof). Then, we will say that $\delta$ is the  {\em height} of $\tcal$,  $\langle |T_{\al}| : \al < \delta \rangle$ is the {\em cardinal sequence} of $\tcal$, and we will write $I_{\al}(\tcal) = T_{\al}$ for $\al < \delta$.

\vspace{1mm}
If $\tcal = \langle T,\preceq,i \rangle$ is an LCS poset and $S\subset T$ such that $i\{s,t\}\subset S$ for all $\{s,t\}\in [S]^2$, we define the {\em restriction of} $\tcal$ to $S$ as $\tcal\upharpoonright S = \langle S, \preceq \upharpoonright (S\times S), \\ i\upharpoonright [S]^2\rangle$.

\vspace{2mm}
\subsection*{ Cardinal sequences and PCF theory} We want to remark that  the topic of cardinal sequences for LCS spaces is related to Shelah's PCF theory.   If $A$ is a set of regular cardinals and $U$ is an ultrafilter over $A$, the set $\prod A/U$ is linearly ordered. Then, we define

$$\mbox{PCF}(A) = \{\mbox{cf}(\prod A/U) : U \mbox{ ultrafilter on } A \}.$$

In  \cite{ru},  Ruyle studied the notion of a PCF structure due to Magidor and Foreman,   a refinement of the notion of an LCS poset,  in which some conditions are added in order to reflect the fundamental properties of the PCF operator on  the set $\{\omega_n : n \geq 1 \}$. The notion of a  PCF structure has interest in the theory of cardinal arithmetic, concerning Shelah's result that if $\aleph_{\omega}$ is a strong limit cardinal then $2^{\aleph_{\omega}} < \aleph_{\omega_4}$ (see \cite{se}). Suppose that  $\tcal = \langle T,\preceq,i \rangle$ is an LCS poset. If $S  \subset T$, we denote by $\overline{S}$ the closure of $S$ in $X(\tcal)$. Then we say that $\tcal$ is a PCF structure,  if   $X(\tcal)$ is a compact space such that the following conditions hold:

\begin{enumerate}

\item $T$ is an infinite successor ordinal.

\item For every $\mu,\xi \in T$, $\mu \prec \xi$ implies $\mu \in \xi$.

\item $\overline{\omega} = T$.

\item If $I\subset T$ is an interval, then $\overline{I}$ is also an interval.

\item For every $\mu\in T$  of uncountable cofinality,  there is a closed unbounded set $C_{\mu}\subset \mu$ such that $\overline{C_{\mu}}\subset \mu + 1$.

\end{enumerate}

   If $\tcal$ is a  PCF structure  and  $X = X(\tcal)$,  we have that $|I_{\alpha}(X)| \leq |\alpha + \omega|$ for every $\alpha < \mbox{ht}(X)$ (see  \cite{ru}). The interest of the notion of  a PCF structure lies in the fact that for the proof of Shelah's theorem stated above, it is shown by means of a combinatorial argument that there is no PCF structure of size $\geq \omega_4$ (see \cite{bm} and \cite{ru}). Then, one could improve Shelah's bound on $2^{\aleph_{\omega}}$ to $\aleph_{\omega_3}$, by showing that in ZFC there is no PCF structure of size $\omega_3$. In \cite{ru}, it was shown that it is consistent with ZFC that there is a PCF structure $\tcal$ such that $X(\tcal)$ has height $\omega_2$, and so we can not hope to improve  Shelah's bound on $2^{\aleph_{\omega}}$ to $\aleph_{\omega_2}$, at least by using PCF theory. This result was improved in \cite{ev}, where it was shown that it is consistent with ZFC  to have PCF structures of height $\delta$ for all ordinals $\delta < \omega_3$. However, it is not known whether it is consistent with ZFC that there is a PCF structure $\tcal$ such that $X(\tcal)$ has height $\omega_3$. If we could answer this question in the affirmative, we could not hope to use PCF theory to improve Shelah's bound on $2^{\aleph_{\omega}}$ to $\aleph_{\omega_3}$. 

However, it is not even known whether it is consistent with ZFC that there is an LCS space $X$ such that $\mbox{CS}(X) = \langle \omega \rangle_{\omega_3}$. This is one of the main open problems in the topic of cardinal sequences.

\section{\bf The forcing construction}

\newcommand{\fnull}{f_0}
\newcommand{\fone}{f_1}
\newcommand{\Znull}{Z_0}
\newcommand{\Zone}{Z_1}

In this section we  carry out the forcing construction that will permit us to prove Theorem 1.3.
Without loss of generality, we may assume that 
$\de$ is a limit ordinal with $\om_2 < \delta < \om_3$.
Fix $f$ as in Theorem~\ref{Thm.1.2}.

\subsection*{Outline of the proof.}

Since the proof is long and complicated, we first sketch the main ideas.

We will define two sequences of cardinals, $\fnull$ and $\fone$,  of length ${\delta}$ such that 
\begin{enumerate}[(F1)]
\item\label{en:FF0F1} $f=\fnull+\fone$, that is, $f({\alpha})=\fnull({\alpha})+\fone({\alpha})$ for each 
${\alpha}<{\delta}$,
\smallskip

\item in some cardinal-preserving generic extension there exist LCS spaces 
$\Znull$ and $\Zone$ with 
$\CS(\Znull)=\fnull$ and $\CS(\Zone)=\fone$.  
\end{enumerate}
Consequently, 
\begin{displaymath}
\CS(\Znull\oplus \Zone)=f, 
\end{displaymath}
and this completes the proof. 

The construction of $\fnull$ is straightforward, and the existence of  $\Znull$ will follow from earlier results.  Define 
\begin{displaymath}
\fnull({\alpha})=\left\{\begin{array}{ll}
{f({\alpha})}&\text{if $\cf({\alpha})<{\omega}_2$,}\\\\
{{\omega}}&\text{if $\cf({\alpha})={\omega}_2$.}\\
\end{array}\right.
\end{displaymath}
Then, by Theorem \ref{tm:8main}, in some cardinal- and cofinality-preserving generic extension 
$V_1$, 
\begin{displaymath}
V_1\models\text{$\fnull=\CS(\Znull)$
 for some LCS space
$\Znull$. 
}
\end{displaymath}

From now on we work in $V_1$. 
The construction of $\fone$ is much more involved. 
Write  
\begin{displaymath}
  L^{\lambda}_{\om_2} = \{\al\in \delta : \cf (\al) = \om_2 \text{ and } \ka_{\al} = \la\}.
\end{displaymath}
For each  ${\beta}\in  L^{\lambda}_{\om_2}$
we will define an ordinal ${\gamma}_{\beta}\le {\beta}$
such that taking 
\begin{displaymath}
  \tilde{L}^{\lambda}_{\om_2}=\{{\beta}\in {L}^{\lambda}_{\om_2}: {\gamma}_{\beta}<{\beta}\} 
\end{displaymath} 
and defining $\fone\in {}^{\delta}\{{\omega},{\lambda}\}$ by the requirement
\begin{equation}\label{eq:f1}
\fone^{-1}\{{\lambda}\}=\bigcup\{[{\gamma}_{\beta},{\beta}]:{\beta}\in \tilde{L}^{\lambda}_{\om_2}\},
\end{equation}
the function $\fone$ will  satisfy condition (F\ref{en:FF0F1}).

\subsection*{Construction of $\fone$}
We may assume that $L^{\lambda}_{\omega_2} \ne \emptyset$, because otherwise $f = \fnull$ and the result is already established.

 Since $f^{-1}\{\om\}$ is $\om_2$-closed in $\delta$,
  if $\al\in L^{\lambda}_{\om_2} $ then there is an ordinal $\al' < \al$ such that for all $\al' \leq \be < \al$ we have 
that $\ka_{\be} = \lambda$.

To define $\fone$ and subsequently carry out our forcing, we will use the notion of a tree of intervals introduced in \cite{ma1}.

 We say that $I$ is an {\em ordinal interval}, if there are ordinals  $\al$ and $\be$  with $\al \leq  \be$ such that $I = [\al,\be) = \{\ga : \al \leq \ga < \be \}$. Then, we write $I^- = \al$ and $I^+ = \be$.

 \begin{definition}  
  \label{Df.2.2}
  Proceeding by induction on $n \in \omega$, we define a family $\mathcal{I}_n$ of ordinal intervals.  
  For $n \ge 1$ and $I \in \mathcal{I}_n$ with $I^{+} \in L^{\lambda}_{\omega_2}$ we also define an ordinal 
  $\gamma(I)$ with $I^- < \gamma(I) \le I^{+}$.
  
  We start with $\mathcal{I}_0 = \{[0,{\delta})\}$.  

  Now suppose $n \ge 0$, and assume that $\mathcal{I}_0, \dots, \mathcal{I}_{n}$ have already been defined and that,
  for every interval $I \in \mathcal{I}_1 \cup \dots \cup \mathcal{I}_{n-1}$ with $I^{+} \in L^{\lambda}_{\omega_2}$, 
  we have defined an ordinal $\gamma(I)$ with $I^- < \gamma(I) \le I^{+}$.
  
  Let $I = [\alpha,\beta) \in \mathcal{I}_n$ with $\beta \in L^{\lambda}_{\omega_2}$ and $n\geq 1$.  
  We define $\gamma(I)$ as follows:
  \begin{enumerate}[$\bullet$]
  \item If there exists an interval $J \in \mathcal{I}_j$ with $0 < j < n$, $J^{+} \in L^{\lambda}_{\omega_2}$,  
  and $I \subset (\gamma(J), J^{+})$, then set $\gamma(I) = \beta$.
  
  \item Otherwise, choose an ordinal $\gamma(I)$ of cofinality $\omega$ with 
  ${\alpha} < \gamma(I) < \beta$ such that  
  $f(\xi) = \lambda$ for all  $\gamma(I) \le \xi \le \beta$.
  \end{enumerate}
    Next, for $I = [\alpha,\beta) \in \mathcal{I}_n$ with $n\geq 0$ we define $E(I)$ and $\intpart{I}$ as follows:
    \begin{enumerate}[(1)]
    
    \item Suppose that $\be$ is a limit ordinal.
    
    \begin{enumerate}[(i)]
  \item If $\beta \in L^{\lambda}_{\omega_2}$, $n > 0$ and $\gamma(I) < \beta$, let  
  $E(I) = \{\epsilon^{I}_{\nu} : \nu < \omega_2\}$
  be a closed cofinal subset of $I$ of order type $\omega_2$ such that  
  \[
  \epsilon^{I}_{0} = \alpha,\qquad 
  \epsilon^{I}_{1} = \gamma(I),\qquad
  \epsilon^{I}_{2} = \gamma(I) + 2,
  \]
  and $\epsilon^{I}_{\nu+1}$ is a successor ordinal for $1 < \nu < \omega_2$.  
  Then set
  \[
  \intpart{I}
     = \{[\epsilon^{I}_{\nu}, \epsilon^{I}_{\nu+1}) : \nu < \omega_2\}.
  \]

  \item Otherwise,
  let
  $E(I) = \{\epsilon^{I}_{\nu} : \nu < \cf(\beta)\}$
  be a closed cofinal subset of $I$ of order type $\cf(\beta)$ with $\epsilon^{I}_0 = \alpha$  
  and $\epsilon^{I}_{\nu+1}$ a successor ordinal for $0 \le \nu < \cf(\beta)$.  
  Then set
  \[
  \intpart{I}
     = \{[\epsilon^{I}_{\nu}, \epsilon^{I}_{\nu+1}) : \nu < \cf(\beta)\}.
  \]
  
   \end{enumerate}
   
  \item Suppose that  $\beta = \beta' + 1$ is a successor ordinal. Then:
  \begin{enumerate}[(i)]
  \item If $I$ is infinite, let $\gamma' < \beta$ be the greatest limit ordinal.
  Define  
$  E(I) = \{\alpha, \gamma', \gamma' + 1, \dots, \beta'\}$
  and
  \[
  \intpart{I}
     = \{[\alpha,\gamma'), \{\gamma'\}, \{\gamma'+1\}, \dots, \{\beta'\}\}.
  \]
  
  \item If $I$ is finite, put  
  $E(I) = \{\alpha,\alpha+1,\dots,\beta'\}$
  and  
  \[
  \intpart{I}
     = \{\{\alpha\},\{\alpha+1\},\dots,\{\beta'\}\}.
  \]
  \end{enumerate}
  \end{enumerate}
  
  Finally, define
  \[
  \mathcal{I}_{n+1}
     = \bigcup \{\intpart{I} : I \in \mathcal{I}_n\}.
  \]
    \end{definition}
  We put  $$\mathbb I=\bigcup\{\ical_n:n<{\omega}\}.$$ Note that $\mathbb
I$ is a cofinal tree of intervals in the sense defined in \cite{ma1}.
So, the following conditions are
satisfied:
\begin{enumerate}[(i)]
\item For every $I,J\in {\mathbb I}$, $I\subset J$ or
$J\subset I$ or $I\cap J = \emptyset$.
\item  If $I,J$ are different elements of ${\mathbb I}$ with
$I\subset J$ and $J^+$ is a limit ordinal, \mbox{ then  $I^+ < J^+$ }.
\item $\ical_n$ partitions $[0,\delta )$ for each $n < \omega$.
\item $\ical_{n + 1}$ refines $\ical_n$ for each $n <\omega$.
\item For every $\alpha <{\delta}$ there is an $I\in {\mathbb I}$
such that $I^- = \alpha$. 
\end{enumerate}

We now make the following definitions:
\begin{enumerate}[(1)]
\item If $\al < \de$ and $n < \om$, we denote by $I(\al,n)$ the interval $I\in  \ical_n$ such that $\al\in I$.
\item If $\al < \de$, we write 
\begin{displaymath}
 \lop({\alpha})=\min\{n\in {\omega}: \exists I\in \mathcal I_n \ (I^-={\alpha})\}.
\end{displaymath}
\end{enumerate}

  The following notion will also be needed in our construction.
 
 \begin{definition}\label{Df.2.3} Assume that $0<\zeta < \delta$ is an ordinal. We define the interval $J(\zeta)$ as follows. 
    If there is an interval $I\in \ical_{\lop(\zeta)}$ with $I^+ = \zeta$, we put $J(\zeta) = I$. 
  Otherwise, $\zeta$ is a limit point of $E(I(\zeta,\lop(\zeta)-1))$, and then we define $J(\zeta) = I(\zeta,\lop(\zeta)-1)$.
\end{definition}
Note that if $\cf (\zeta) = \om_2$, then
 ${\zeta}\in E(I(\zeta,\lop(\zeta)-1))$
but $|E(I(\zeta,\lop(\zeta)-1))\cap \zeta|\le \om_1$, hence $\zeta$ can not be a limit point of $E(I(\zeta,\lop(\zeta)-1))$. Thus, the ordinal $\zeta$ has an immediate predecessor $\xi$ in $E(I(\zeta,\lop(\zeta)-1))$, and so we have $J(\zeta) = [\xi,\zeta) \in \ical_{\lop(\zeta)}$. Analogously, if $\zeta$ is a successor ordinal, it is clear that $J(\zeta) = [\xi,\zeta)$ where $\xi$ is the immediate predecessor of  $\zeta$ in $E(I(\zeta,\lop(\zeta)-1))$.

 Assume that $\be\in L^{\lambda}_{\om_2}$. Note that if  $[\al,\beta') = I(\be,\lop(\be)-1)$ then, by Definition 2.1, $\be'$ is a successor ordinal $\zeta +1$ and  $E([\al,\beta')) = \{\al,\be, \be + 1\ ,\dots, \zeta \}$. Then, we define 
 \begin{equation}\notag
  \ga_{\be} = \ga(J(\be)) = \ga([\al,\be)).
 \end{equation}
 Recall that
 $\tilde{L}^{\lambda}_{\om_2} = \{\be \in L^{\lambda}_{\om_2} : \ga_{\be} < \be \}$.
 { Thus we have defined $\fone$  by \eqref{eq:f1}.  }
 Since $\bigcup\{[{\gamma}_{\beta},{\beta}]:{\beta}\in \tilde{L}^{\lambda}_{\om_2}\} \supset
 {L}^{\lambda}_{\om_2}$, $\fone$ satisfies requirement (F1).

\vspace{2mm}
\subsection*{The forcing construction of $\Zone$}
The remainder of the paper will be  devoted to construct a c.c.c  generic extension of $V_1$ which contains an LCS space $\Zone$
with $\CS(\Zone)=\fone.$ 

Since $\square_{\om_1}$ is preserved by any cardinal-preserving forcing, 
we have a $\Delta$-function in $V_1$ by Theorem \ref{Thm.1.4}.

The following basic conditions should be noted:

\begin{enumerate}[\bigskip\hspace{4.5mm} ($1$)]
\item[($*$) (1)]    For every $\be\in L^{\lambda}_{\om_2}\setminus \tilde{L}^{\lambda}_{\om_2}$ there is $\be'\in  \tilde{L}^{\lambda}_{\om_2}$ such that $\be \in (\ga_{\be'},\be')$.
\smallskip\addtocounter{enumi}{1}

\item If $\be\in \tilde{L}^{\lambda}_{\om_2}$, then $\lop(\ga_{\be}) = \lop(\be) + 1$, $I(\be,\lop(\be)) = \{\be\}$ and \\ $I(\ga_{\be},\lop(\be) + 1) = \{\ga_{\be}, \ga_{\be} + 1\}$. \smallskip

\item  If $\be,\be'\in \tilde{L}^{\lambda}_{\om_2}$ with $\be \neq \be'$, then $[\ga_{\be},\be]\cap [\ga_{\be'},\be']  = \emptyset$. \smallskip

\item If $\be\in \tilde{L}^{\lambda}_{\om_2}$ and $J(\beta)\in \ical_n$, then $J(\be) = I(\ga_{\be},n)$,  $J(\ga_{\be})\in \ical_{n+1}$ and $[\ga_{\be}, \ga_{\be} + 2)\in \ical_{n+1}$.
\smallskip

\item If $\be\in \tilde{L}^{\lambda}_{\om_2}$, then $J(\beta) = [ I(\be, l(\be) - 1)^-,\be)$.
\smallskip

\item    If $\be\in \tilde{L}^{\lambda}_{\om_2}$, then $I(\ga_{\be}, l(\be)-1)  = I(\be, l(\be)-1)$  and $J(\ga_{\be}) = [ I(\be, l(\be) - 1)^-,\ga_{\be})$.
\smallskip

\item If $\be\in \tilde{L}^{\lambda}_{\om_2}$ and $I(\ga_{\be}, n)\supsetneqq  J(\beta)$, then $I(\ga_{\be}, n+1) \supsetneqq [\ga_{\be},\ga_{\be} + 2)$.
\smallskip

\item If $\be\in \tilde{L}^{\lambda}_{\om_2}$,  then  for every $I \in \mathbb I$ with $[\ga_{\be}, \ga_{\be} + 2)\subsetneqq I$ we have that $\ga_{\be} > I^-$.

\end{enumerate}
Conditions $(*)(1)-(4)$ are clear.  To check $(*)(5)$, put $\al =  I(\be, l(\be)-1)^-$. Then, the  first element of $E(I(\be, l(\be)-1))$ is $\al$ and its second element is $\be$, and so $J(\beta) = [ I(\be, l(\be) - 1)^-,\be)$.

\vspace{1mm}
To verify $(*)(6)$, note that $\ga_{\be}\in I(\be,i)$ for every $i < l(\beta)$, and so $I(\ga_{\be},l(\be) - 1) = I(\be,l(\be) - 1)$. Put $\al =  I(\be, l(\be)-1)^-$.  Then, in $E(J(\beta))$ the first element is $\al$ and the second element is $\ga_{\be}$, and hence $J(\ga_{\be}) = [\al,\ga_{\be})$.

\vspace{1mm} To verify $(*)(7)$, note that if $I(\ga_{\be}, n + 1) = \{\ga_{\be}\}$ then $I(\ga_{\be}, n) =  [\ga_{\be}, \ga_{\be} + 2)$, and if
$I(\ga_{\be}, n + 1) = [\ga_{\be}, \ga_{\be} + 2)$ then  $I(\ga_{\be}, n) = J(\beta)$.

\vspace{1mm} And to check $(*)(8)$, note that 
if $[\ga_{\be}, \ga_{\be} + 2)\subsetneqq I$, then $J(\be)\subset I$ by $(*)(4)$, and so  $I^- \leq J(\be)^- < \ga_{\be}$.

\begin{definition}\label{df:k}

Assume that $\al < \be < \de$ and $I\in  {\mathbb I}$. We say that $\al$ and $\be$ {\em separate at } $I$ if, for some $n < \om$, we have $I = I(\al,n) = I(\be,n)$ but $I(\al, n+1) \neq I(\be,n+1)$. In this case we say that $n$ is the {\em level where} $\al$ and $\be$ {separate}, and we write $\kop(\al,\be) = n$, i.e. 
\begin{displaymath}
\kop(\al,\be)=\min\{n\in {\omega}: I(\al, n+1) \neq I(\be,n+1)\}.
\end{displaymath}
 
\end{definition}
Clearly, any pair of different elements of $\de$ separate at some interval of ${\mathbb I}$.

\begin{definition} Assume that $\al$ and $\be$ are ordinals with $\al < \be < \de$. Let $k = \kop(\al,\be), J = I(\al, k+1)$ and $K = I(\be,k+1)$. For $n <\om$, we put $\al_n = I(\al,n)^+$ and $\be_n = I(\be,n)^-$. Then, we define  the {\em walk} from $\al$ to $\be$ via ${\mathbb I}$, denoted $w(\al,\be)$, as follows. 

If $\lop(\al) \leq k + 1$, then 
$$ w(\al,\be) = \langle \al,  \al_{k+ 1}, \be_{k+1},\be_{k+2}, \dots, \be_{\lop(\be) - 1},\be \rangle.$$

If $\lop(\al) > k + 1$, then 
$$ w(\al,\be) = \langle \al, \al_{\lop(\al)}, \al_{\lop(\al) -1},\dots, \al_{k+1}, \be_{k+1},\be_{k+2}, \dots, \be_{\lop(\be) -1} ,\be \rangle.$$
\end{definition}

Observe that the elements of a walk are not necessarily pairwise distinct.

 This notion of a walk is a refinement of the notion of a walk introduced in \cite[Section 2]{ma1}, and in a sense is similar to the notion of a walk introduced by Todor\v{c}evi\'c in \cite{to1}.

\smallskip
The following fact should be noted.

\begin{proposition}\label{Pr.2.5} If $\be \in \tilde{L}^{\lambda}_{\om_2}$ and $\al < \ga_{\be}$, then $\ga_{\be}$ is in the sequence  $w(\al,\ga_{\be} + 1)$.
 \end{proposition}

\vspace{2mm}\noindent {\bf Proof}.
Let $k = \kop(\al,\ga_{\be} + 1)$ and let $I$ be the interval where $\al$ and $\ga_{\be} + 1$ separate.  We may assume that $\lop(\al) > k + 1$. Otherwise, the considerations are simpler.  Let $J = I(\al, k+1)$ and $K = I(\ga_{\be} + 1,k+1)$.  Note that $I \supset J(\be)$, because if $J(\be)\in \ical_n$ we see that $I(\ga_{\be}, n +1) = [\ga_{\be},\ga_{\be} + 2)$, and so $k\leq n$. First, assume that $I = J(\be)$. Then, $J=J(\ga_{\be})$ and $K =[\ga_{\be}, \ga_{\be} + 2)$, and so $J^+ = K^- = \ga_{\be}$. Then, $w(\al,\ga_{\be} + 1) = \langle \al,  \al_{\lop(\al)}, \al_{\lop(\al) -1},\dots, \ga_{\be}, \ga_{\be} + 1\rangle$.

 Now, suppose that $I \supsetneqq J(\be)$. So,  $\al < J(\be)^-$.  For $n <\om$, we write $\al_n = I(\al,n)^+$ and $\be_n = I(\be,n)^-$. Then, $ w(\al,\be) = \langle \al, \al_{\lop(\al)}, \al_{\lop(\al) -1},\dots, \al_{k+1},$ $\be_{k+1},\be_{k+2}, \dots, \be_{\lop(\be) -1}, \ga_{\be}, \ga_{\be}  + 1 \rangle$, where $\be_{\lop(\be) -1} = J(\be)^-$ by condition $(*)(5)$.  $\square$

   However, note that if  $\be \in \tilde{L}^{\lambda}_{\om_2}$  and $\al < J(\be)^-$, then $\ga_{\be} $ does not appear  in the sequence  $w(\al,\be)$, because $\lop(\ga_{\be}) > \lop(\be)$.
  
  \begin{proposition}\label{Pr.2.6}
     Assume that $\be\in \tilde{L}^{\lambda}_{\om_2}$. Let $I = J(\ga_{\be})$ and  $E(I) = \langle \epsilon_n : n < \om \rangle$. Let $\al\in I$ and let $m$ be the least natural number $n$ such that $\al < \epsilon_n$. Then:
  \begin{enumerate}[(a)]
  \item $\epsilon_m$ is in the sequence $w(\al,\ga_{\be})$.
  \item $\epsilon_m$ is in the sequence $w(\al,\ga_{\be} + 1)$.
  \end{enumerate}
  \end{proposition}

  \vspace{2mm}\noindent {\bf Proof}.
 We prove condition (b); the proof of (a) is similar.  Let $k = \kop(\al, \ga_{\be} + 1)$. Note that $J(\be)$ is the interval where $\al$ and $\ga_{\be} + 1$ separate by condition $(*)(4)$. Hence, $I(\al,k+1) = J(\ga_{\be})$ and $I(\ga_{\be} + 1,k+1) =  \{\ga_{\be},\ga_{\be} + 1\}$. First, assume that $\al > \epsilon_{m-1}$. 
 Then,
 $w(\al,\ga_{\be} + 1) = \langle \al, \al_{l(\al)}, \al_{l(\al) -1},\dots, \epsilon_m,\ga_{\be},\ga_{\be} + 1\rangle$.  And if $\al = \epsilon_{m-1}$, then $w(\al,\ga_{\be} + 1) = \langle \al, \epsilon_m, \ga_{\be},\ga_{\be} + 1\rangle$.   $\square$

\begin{definition}
  For every $n \in \omega$ and every $I \in \mathcal{I}_n$, we choose a function 
  \begin{equation}\label{eq:GI}\notag
    G_I : [E(I)]^2 \longrightarrow [E(I)]^{\le \omega}  
  \end{equation}
  as follows.
 If $\cf(I^+) \le \omega_1$, then for every $\{\alpha,\beta\} \in [E(I)]^2$ we set
$  G_I\{\alpha,\beta\} = E(I)\cap {\alpha}\cap {\beta}.$ 
 If $\cf(I^+) = \omega_2$, then $G_I$ is a $\Delta$-function on $E(I)$, which exists  
  by Theorem~\ref{Thm.1.4}.
  Next, we define the function 
  \[
  G : [\delta]^2 \longrightarrow [\delta]^{\le \omega}
  \]
  as follows.  
    If $\alpha < \beta < \delta$, let $I \in \mathcal{I}_n$ be the interval where $\alpha$ and $\beta$ separate, 
  and let 
$  J = I(\alpha, n+1)$ and $K = I(\beta, n+1)$.
  Then we set
  \[
  G\{\alpha,\beta\} = G_I\{J^-, K^-\} \cup \{I^-\}.
  \]
  \end{definition}

\medskip

In order to define the required notion of forcing, we need some additional preparation.

The desired space $\Zone$ will be constructed in two steps. 
First, we add an LCS poset $\mathcal T$ to our model, and then, 
 in Lemma \ref{lm:Zone}
we show how $\Zone$ can be obtained from the LCS poset $\mathcal T$. 
We start by defining the set $Y$, which will serve as the underlying set of the LCS poset  $\mathcal T$.
The set $Y$ will be obtained as a union of blocks.
First, we define the index set of our blocks as
\begin{displaymath}
  {\mathbb B }  = \{1\} \cup \{\langle \al,\zeta \rangle : \al\in \tilde{L}^{\lambda}_{\om_2} \mbox{ and } 0 < \zeta < \lambda \}.
\end{displaymath}
Let  $${\mathbb B }_1 = \delta \times \om,$$ and for $\al\in \tilde{L}^{\lambda}_{\om_2}$ and $0 < \zeta < \lambda$ define
$${\mathbb B }_{\langle \al,\zeta \rangle} = (\{\ga_{\al}\}\times [\om\cdot \zeta, \om\cdot \zeta + \om)) \cup \{\langle \ga_{\al} + 1 ,\lambda + \zeta \rangle \}.$$
Then, we define 
\begin{displaymath}
Y=\mathbb{B}_1\cup\bigcup\{\mathbb B_{\<{\alpha},{\zeta}\>}:{\alpha}\in \tilde{L}^{\lambda}_{\om_2}
\text{ and }0<{\zeta}<{\lambda}\}.
\end{displaymath}
 
We define the function $$\pi_{{\mathbb B}} : Y \rightarrow {\mathbb B}$$ by setting $\pi_{\mathbb{B}}(x) = B$ if and only if $x \in \mathbb{B}_B$, for $x \in Y$.

We also define the functions
\begin{displaymath}
{\pi}:Y\to {\delta}\text{ and }{\pi}_-:Y\to {\delta}
\end{displaymath}
as follows: for $s = \langle \alpha, \eta \rangle \in Y$, set $\pi(s) = \alpha$.  

For $\pi_-$, set $\pi_-(s) = \pi(s)$ if $s \in \mathbb{B}_1$, and $\pi_-(s) = \gamma_{\alpha}$ if $s \in \mathbb{B}_{\langle \alpha, \zeta \rangle}$.

\smallskip

 We write
$$U = \bigcup \{ {\mathbb B }_{\langle \al,\zeta \rangle} : \al\in \tilde{L}^{\lambda}_{\om_2} \mbox{ and } 0 < \zeta < \lambda \}=Y\setm \mathbb B_1.$$

\vspace{1mm}
  The following definition will be essential in our notion of forcing.

\begin{definition}
 We define the function $h : [Y]^2 \rightarrow [\de]^{\leq \om}$ as follows. If $\{s,t\}\in [Y]^2$ with $\pi_{-}(s) \leq \pi_{-}(t)$  we put

\begin{equation}\notag
 h\{s,t\} =\left\{
 \begin{array}{ll}
 G\{\pi_{-}(s),\pi_{-}(t)\}  &\text{if  $ \pi_{-}(s) < \pi_{-}(t)$,}\medskip\\
 E(J(\pi_{-}(s)))  &\text{if  $\pi_{-}(s) = \pi_{-}(t)$ and
  $\pi_{\mathbb B}(s) \neq \pi_{\mathbb B}(t)$,}\medskip\\
  \emptyset &\text{if  $\pi_{-}(s) = \pi_{-}(t)$ and
  $\pi_{\mathbb B}(s) = \pi_{\mathbb B}(t)$.}
 \end{array}
 \right.
\end{equation}
\end{definition}

\vspace{2mm}
Now, by using the function $h$, we define the desired  notion of forcing. We write $U_1 = \{s\in U : \pi(s) = \pi_{-}(s) \}$  and $U_2 = \{s\in U : \pi(s) \neq \pi_{-}(s) \}$.

\begin{definition} 
    We define $P_{\delta}$ as the set of all $p = \langle x_p, \preceq_p,i_p \rangle$ such that the following conditions hold:

\begin{enumerate}[(P1)]
\vspace{1mm}
\item $x_p$ is a finite subset of $Y$,

\vspace{1mm}
\item $\preceq_p$ is a partial order on $x_p$ such that:

\begin{enumerate}[(a)]
\item $s \prec_p t$ implies $\pi(s) < \pi(t)$,
\item if $s\prec_p t$ and $s\in U$ then $\pi_{{\mathbb B}}(s) = \pi_{{\mathbb B}}(t)$,
\item if $s \prec_p t$, $s\in {\mathbb B}_1$ and $t\in U_2$, then there is $u\in x_p \cap U_1$ with  $\pi_{\mathbb B}(u) = \pi_{\mathbb B}(t)$ such that $s \prec_p u \prec_p t$,
\item if $x_p\cap {\mathbb B}_{\langle \al,\zeta\rangle}\neq \emptyset$ for $\al \in \tilde{L}^{\lambda}_{\om_2}$ and $0 < \zeta < \lambda$, then $\langle \ga_{\al} + 1, \lambda + \zeta\rangle \in x_p$ and $s \preceq_p \langle \ga_{\al} + 1, \lambda + \zeta\rangle$ for all $s\in x_p\cap {\mathbb B}_{\langle \al,\zeta\rangle}$.

\end{enumerate}

\vspace{1mm}
\item $i_p : [x_p]^2 \rightarrow [x_p]^{<\om}$ satisfies the following:

\begin{enumerate}[(a)]
\item if $s \prec_p t$ then $i_p\{s,t\} = \{s\}$,
\item if $s,t\in x_p\cap {\mathbb B}_1 $ with $s\neq t$ and $\pi(s) = \pi(t)$  then  $i_p\{s,t\} = \emptyset$,
\item if $s,t\in x_p\cap U_1$ with $s\neq t$ and $\pi_{{\mathbb B}}(s) =\pi_{{\mathbb B}}(t)$ then $i_p\{s,t\} = \emptyset$,
\item if $s,t\in x_p$ are $\preceq_p$-incomparable but $\preceq_p$-compatible then $\pi[i_p\{s,t\}] \subset h\{s,t\}$,
\item $v\preceq_p s,t$ for all $v\in i_p\{s,t\}$,
\item for every $u\preceq_p s,t$ there is $v\in i_p\{s,t\}$ such that $u\preceq_p v$.

\end{enumerate}

\vspace{1mm}
\item  If $s <_p t$ with $s\in {\mathbb B}_1$, then there are $u_1,\dots , u_n \in x_p$ with $s \prec_p u_1 \preceq_p \dots \preceq_p  u_n \preceq_p t$ such that
$$w(\pi(s),\pi(t)) = \langle \pi(s),\pi(u_1),\dots , \pi(u_n),\pi(t)\rangle.$$

\end{enumerate}
Now, if $p,q\in P_{\de}$, we write $p\leq_{\de} q$ iff $x_p\supset x_q$, $\preceq_q=\preceq_p \cap (x_q\times x_q)$ and
$i_p\supset i_q$.

We put ${\mathbb P}_{\delta} = \langle P_{\de},\leq_{\de} \rangle$.
\end{definition}

 We will frequently use condition $(P2)(b)$ without explicit mention.

\smallskip Suppose that $p\in   P_{\de}$ and $s \prec_p t$ with $s\in {\mathbb B }_1 $. Then,  if $u_1,\dots , u_n$ are elements of $x_p$ satisfying condition $(P4)$ for $s,t$, we say that the sequence $\langle s, u_1, \dots , u_n, t\rangle$ is a ``{\em walk from $s$ to $t$ in $p$}''. Note that if  $\langle s_1,\dots , s_n\rangle$ is a walk from $s_1$ to $s_n$ in $p\in   P_{\de}$, we may have $s_i = s_{i+1}$ for some $i\in \{1,\dots,n-1 \}$. Then, if $\langle s_1,\dots, s_i,u,\dots ,u, s_j, \dots , s_n\rangle$ is a walk from $s_1$  to $s_n$ in $p$, we will also consider the sequence $\langle s_1,\dots, s_i,u, s_j, \dots , s_n\rangle$ as a walk from $s_1$  to $s_n$ in $p$.

\medskip

Observe that the following  lemmas clearly complete the proof of Theorem \ref{Thm.1.2}.

\begin{lemma}\label{Lm.2.10} ${\mathbb P}_{\delta}$ is c.c.c.

\end{lemma}

\begin{lemma}\label{Lm.2.11} Forcing with  ${\mathbb P}_{\delta}$ adjoins an LCS poset on $Y$.

\end{lemma}

\begin{lemma}\label{lm:Zone}
  Forcing with  ${\mathbb P}_{\delta}$ adjoins an LCS space $\Zone$ with 
  cardinal sequence $\fone.$
\end{lemma}

We will prove first Lemma  \ref{lm:Zone}.
In order to prove Lemmas \ref{Lm.2.10} and \ref{Lm.2.11}, we will refine the arguments given in \cite{ma1}.  
The proof of Lemma~ \ref{Lm.2.10} is particularly intricate and will therefore be presented in full in Section \ref{sc:Lm.2.10}.

\vspace{2mm}\noindent {\bf Proof  of Lemma 2.12.}
Consider a ${\mathbb P}_{\de}$-generic extension $V_1[G]$. 
Clearly, as $|P_{\de}| = \la$ and ${\mathbb P}_{\de}$ is c.c.c., we have that $2^{\om} = \lambda$ in $V_1[G]$. 
From now on we work in $V_1[G]$. Let  $\tcal = \langle Y,\preceq,i \rangle$ where  
$\preceq = \bigcup \{\preceq_p : p\in G \}$ and $i = \bigcup \{i_p : p\in G \}$. 
By Lemma 2.11, $\tcal$ is an LCS poset on $Y$.  Our purpose is to construct from   $\tcal$ an LCS poset 
$\tcal'$ whose underlying set is $(\delta \times \om ) \cup \bigcup \{[\ga_{\al},\al] \times (\la\setminus \om) : \al \in \tilde{L}^{\lambda}_{\om_2} \}$ 
 in such a way that 
 if $\Zone$ is the LCS space associated with $\tcal'$, then 
 $\CS(\Zone)=\fone$.

         In order to construct $\tcal'$,  we  will replace the non-standard blocks ${\mathbb B }_{\langle \al,\xi  \rangle}$ with new blocks ${\mathbb B }'_{\langle \al,\xi  \rangle}$, which are defined as follows. If $\al\in \tilde{L}^{\lambda}_{\om_2}$ and $0 < \xi < \lambda$ we put
       $${\mathbb B }'_{\langle \al,\xi  \rangle} = C^{\xi}_{\al} \cup  S^{\xi}_{\al},$$
with  $C^{\xi}_{\al} = \{\gamma_{\al}\}\times [\om\cdot \xi,\om\cdot \xi + \om)$ and $S^{\xi}_{\al} = ([\ga_{\al} + 1,\al) \times [\om\cdot \xi, \om\cdot \xi + \om ))\cup \{ t^{\al}_{\xi}\} $ where $t^{\al}_{n + 1} = \langle \al, \om + n \rangle$ for $n\in \om$ and $t^{\al}_{\xi} = \langle \al, \om + \xi \rangle$ for $\om\leq \xi < \lambda$.

       If $\al\in \tilde{L}^{\lambda}_{\om_2}$, we put $\delta_{\al} = \mbox{o.t.}(\al\setminus \gamma_{\al}) + 1$. 
 Note that since $\tcal \upharpoonright  {\mathbb B}_1$ is an LCS poset whose cardinal sequence is $\langle \om \rangle_{\de}$, for every $\al\in \tilde{L}^{\lambda}_{\om_2}$ and $0 < \xi < \lambda$ we can define an LCS poset $\tcal^{\xi}_{\al} = \langle {\mathbb B}'_{\langle \al,\xi  \rangle}, \preceq^{\xi}_{\al} , i^{\xi}_{\al}\rangle$ of height $\delta_{\al} + 1$ whose cardinal sequence is $\langle \om \rangle_{\delta_{\al}}$ such that  $I_0(\tcal^{\xi}_{\al}) = C^{\xi}_{\al}$ and $I_{\delta_{\al}}(\tcal^{\xi}_{\al}) = \{t^{\al}_{\xi} \}$. Then, we define the required poset $\tcal' =  \langle Y', \preceq',i' \rangle$ as follows. We put

      $$Y' = {\mathbb B }_1  \cup \bigcup \{{\mathbb B }'_{\langle \al,\xi  \rangle}  : \al \in \tilde{L}^{\lambda}_{\om_2} \mbox{ and } 0 < \xi < \lambda \}.$$

      \noindent For every $s,t\in Y'$, we put $s \prec' t$ iff one of the following conditions holds:

      \begin{enumerate}[(a)]
      \item $s,t\in Y$ and $s \prec t$,
      \item $s,t\in {\mathbb B }'_{\langle \al,\xi  \rangle}$ for some $\al \in \tilde{L}^{\lambda}_{\om_2}$ and $0 < \xi < \lambda$ and $s \prec^{\xi}_{\al} t$,
      \item $s\in {\mathbb B}_1, t\in {\mathbb B }'_{\langle \al,\xi  \rangle} $ for some  $\al \in \tilde{L}^{\lambda}_{\om_2}$ and $0 < \xi < \lambda$ and there is $u\in C^{\xi}_{\al}$ such that $s \prec u \preceq^{\xi}_{\al} t$.

      \end{enumerate}

      \noindent It is easy to see that $\preceq'$ is a partial order on $Y'$. Now, we define the infimum function $i'$. Assume that $\{s,t\}\in [Y']^2$. If $s \prec' t$, we put $i'\{s,t\} = \{s\}$, and if $t \prec' s$, we put $i'\{s,t\} = \{t\}$. So, suppose that $s,t$ are $\preceq'$-incomparable. We distinguish the following cases.

\newcases
\begin{case}
  $s,t\in Y$.
\end{case}
We put $i'\{s,t\} = i\{s,t\}$.

      \begin{case}
        $s,t\in  {\mathbb B }'_{\langle \al,\xi  \rangle}$ for some $\al\in \tilde{L}^{\lambda}_{\om_2}$ and $0 < \xi < \lambda$.
      \end{case}

      We put $i'\{s,t\} = i^{\xi}_{\al}\{s,t\}$. 
      
      \vspace{1mm}
       It is easy to check that $i'\{s,t\}$ is correct in Cases 1 and 2, i.e. it satisfies condition $(4)$ of Definition 1.6.

\begin{case}
  $s\in  S^{\xi}_{\al}$ and $t\in  S^{\eta}_{\be}$  for some $\al,\be\in \tilde{L}^{\lambda}_{\om_2}$ and $0 < \xi,\eta < \lambda$ in such a way that either ($\al \neq \be$) or ($\al = \be$ and $\xi\neq\eta$).
\end{case}

       Let $s^* = \langle \ga_{\al} + 1, \lambda + \xi \rangle$ and $t^* = \langle \ga_{\be} + 1, \lambda + \eta \rangle$. So, $s^*$ and $t^*$ are the top points of the blocks ${\mathbb B }_{\langle \al,\xi  \rangle}$ and ${\mathbb B }_{\langle \beta,\eta  \rangle}$  respectively. By condition $(P2)(b)$, $s^*$ and $t^*$ are $\preceq$-incomparable. Then, we define

      $$i'\{s,t\} = \{v\in i\{s^*,t^*\}:  v \prec'  s \mbox{ and }  v \prec'  t \}.$$
      
      \vspace{1mm} We show that $i'\{s,t\}$ is correct. So, suppose that $v \prec' s,t$. Clearly, $v\in {\mathbb B}_1$.
       So, there are $u_1\in C^{\xi}_{\al}$ such that $v \prec u_1 \prec^{\xi}_{\al} s$ and $u_2\in C^{\eta}_{\be}$ such that $v \prec u_2 \prec^{\eta}_{\be} t$. Hence, $v\prec s^*,t^*$. So, there is $w\in i\{s^*,t^*\}$ such that $v \preceq w$. Now, by using condition $(P2)(c)$, there are $u'_1\in C^{\xi}_{\al}$ and $u'_2\in C^{\eta}_{\be}$ such that $w \prec u'_1 \prec   s^*$ and $w \prec u'_2 \prec   t^*$. Now, since $v \preceq w$, by using condition $(P3)(c)$, we infer that $u_1 = u'_1$ and $u_2 = u'_2$. Thus
      $w\prec u_1$ and $w\prec u_2$, and therefore $w\in i'\{s,t\}$.

\begin{case}
  $t\in  S^{\xi}_{\al}$ for some $\al\in \tilde{L}^{\lambda}_{\om_2}$ and $0 < \xi < \lambda$ and $s\in Y\setminus  C^{\xi}_{\al}$.
\end{case}

Let $t^* = \langle \ga_{\al} + 1, \lambda + \xi \rangle$.  By condition $(P2)(b)$, $t^*\not\preceq s$. Then, if  $s \prec t^*$  we define $i'\{s,t\} = \emptyset$. Otherwise, we put
      $$i'\{s,t\} = \{v\in i\{s,t^*\}:  v \prec'  t \}.$$
       Again, we show that $i'\{s,t\}$ is correct. If $s$ and $t^*$ are $\preceq$-incomparable,  we can use an argument similar to the one given in Case 3. So, suppose that $s \prec t^*$. Then, we show that there is no $v\in Y'$ such that $v\prec' s,t$. For this, assume on the contrary that there is $v\in Y'$ with  $v\prec' s,t$.  Since $v \prec' t$, there is $u\in C^{\xi}_{\al}$ such that $v \prec u \prec^{\xi}_{\al} t$. But as $s \prec t^*$, there is $u'\in  C^{\xi}_{\al}$ such that $s \prec u'$ by $(P2)(c)$. So as $v\prec u,u'$, we infer that $u=u'$ by $(P3)(c)$. Thus, we deduce that $s \prec' t$, which contradicts the assumption that $s$ and $t$ are $\preceq'$-incomparable.

         \vspace{1mm}

         Now, it is easy to check that $\tcal'$ is as required. This completes the proof of Lemma \ref{lm:Zone}.

 \vspace{2mm}
 \noindent
{\bf Proof of Lemma 2.11.}
 Let $G$ be a ${\mathbb P}_{\delta}$-generic filter. We write $\preceq = \bigcup \{\preceq_p : p\in G \}$ and
$i = \bigcup \{i_p : p\in G \}$.  First, for every $y\in Y$, we show that
$\{p\in P_{\delta} : y\in x_p\}$ is dense in ${\mathbb P}_{\delta}$, and hence $Y = \bigcup \{x_p : p\in G \}$. So, assume that $y\in Y$ and $p\in P_{\de}$. We suppose that
 $y\not\in x_p$ and $y\in U_1$. Otherwise, the argument is similar.
  So, let $\al \in \tilde{L}^{\lambda}_{\om_2}$ and $0 < \zeta < \lambda$ be such that $y\in {\mathbb B}_{\langle \al,\zeta \rangle}$. Then, we define the condition $q = \langle x_q,\preceq_q,i_q\rangle$  by $x_q = x_p \cup \{y,\langle \ga_{\al} + 1,\lambda + \zeta\rangle\}$, $\prec_q = \prec_p \cup \{\langle y,\langle \ga_{\al} + 1,\lambda + \zeta\rangle \rangle\}$ and $i_q\{s,t\} = i_p\{s,t\}$ if $\{s,t\}\in [x_p]^2$, $i_q\{s,t\} = y$ if $s = y$ and $t = \langle \ga_{\al} + 1,\lambda + \zeta\rangle $ and  $i_q\{s,t\} = \emptyset$ otherwise. Clearly, $q \leq_{\de} p$.

\vspace{1mm} Now, we prove that $\langle Y, \preceq , i \rangle$ satisfies condition $(3)$ in Definition 1.6. If $\al < \be < \delta$, $t\in Y$ with $\pi(t) = \be$ and $j< \om$ let

 $$D_{t,\al,j} = \{q\in P_{\de}: t \in x_q \mbox{ and for some } n >  j, \langle \al, n \rangle \prec_q t \}.$$

 \noindent We prove that $D_{t,\al,j}$ is dense in  ${\mathbb P}_{\delta}$. So, assume that $\al < \be < \de$, $t\in Y$ with $\pi(t) = \be$, $j < \om$ and $p\in P_{\de}$.
  Our aim is to show that there is $q\in D_{t,\al,j}$ with $q  \leq_{\de} p$. We suppose that $t\in x_p \cap {\mathbb B}_1$. Otherwise, the argument is simpler, because if $t\in U$ then $|\{v\in x_p : t \prec_p v \}| \leq 1$. Let $n$ be a natural number with $n > j$  such that for every $u = \langle \gamma, l \rangle\in x_p \cap {\mathbb B}_1$, $n > l$. We put $s = \langle \al,n \rangle$. Let $k = \kop(\al,\be)$. We  assume that $\lop(\al) > k+1$. Otherwise, the argument is simpler.
  For $k<l\leq \lop(\al)$  we write $s_l = \langle \al_l,n\rangle$ where $\al_l = I(\al,l)^+$. Also, for $k<l<\lop(\be)$ we write $t_l = \langle \be_l,n \rangle$ where $\be_l = I(\be,l)^-$. Then, we define $q = \langle x_q,\preceq_q,i_q \rangle$ as follows:

\begin{enumerate}[(1)]
\item $x_q = x_p \cup \{s\} \cup \{s_l : k < l \leq \lop(\al) \} \cup \{t_l : k < l < \lop(\be) \}$.
\vspace{1mm}
\item $\prec_q = \prec_p \cup$ $\{\langle s,v\rangle : t\leq_p v  \} \cup \{\langle s_l,v \rangle : k < l \leq \lop(\al), t\preceq_p v \}  \cup \{\langle t_l,v \rangle : k < l < \lop(\be), t\preceq_p v \}  \cup \{\langle s, s_l \rangle : k < l \leq \lop(\al)\}  \cup \{\langle s,t_l \rangle :$ $ k < l < \lop(\be) \}  \cup \{\langle s_l,s_m \rangle : m < l \leq \lop(\al) \}  \cup \{\langle t_l,t_m \rangle : l < m < \lop(\be) \}  \cup \{\langle s_l,t_m \rangle : k < l \leq \lop(\al), k < m < \lop(\be) \}$.

  \vspace{1mm}
    \item $i_q\{u,v\} = i_p\{u,v\}$ if $\{u,v\}\in [x_p]^2$;  $i_q\{u,v\} = \{u\}$ if $u \prec_q v$; $i_q\{u,v\} = \{v\}$ if $v \prec_q u$; and  $i_q\{u,v\} = \emptyset$ otherwise.
    \end{enumerate}

    \vspace{2mm} Our purpose is to show that $q\in P_{\de}$. The verification of $(P1),(P2)$ and $(P3)$ is easy. We prove condition $(P4)$. Suppose that $u \prec_q v$ with $u\in  {\mathbb B }_1$. Note that if $\{u,v\}\in [x_p]^2$, we are done because $p$ satisfies $(P4)$. So, assume that $\{u,v\}\not\in [x_p]^2$   First, suppose that $t \not\preceq_p v$. Then, if $u = s_l$ and $v = s_m$ for $m < l$, we have that $\langle s_l,s_{l-1},\dots , s_m \rangle$ is a walk from $s_l$ to $s_m$ in $q$; if   $u = t_l$ and $v = t_m$ for $l < m$, then $\langle t_l,t_{l+1},\dots , t_m \rangle$ is a walk from $t_l$ to $t_m$ in $q$; and if $u = s_l$ and $v= t_m$ for $l \leq \lop(\al)$ and $m < \lop(\be)$, we have that $\langle s_l,s_{l-1},\dots , s_{k+1},t_{k+1},t_{k+2}, \dots, t_m \rangle$ is a walk from $s_l$ to $t_m$ in $q$. And we proceed in a similar way if $u = s$ and $v\in \{s_l : k < l \leq \lop(\al) \} \cup  \{t_l : k < l < \lop(\be) \}$.

     \vspace{1mm} Now, assume that $t\preceq_p v$. We may assume that $u = s$ and $t \prec_p v$. Otherwise, the considerations are analogous. We put $\zeta = \pi(v)$ and $m= \kop(\al,\zeta)$. Note that $m\leq k$. For this, assume on the contrary that $m > k$. Then, as $m = \kop(\al,\zeta)$ and $k+1\leq m$, we infer that $\al,\be,\zeta \in I(\al, k+1)$, which contradicts the fact that $k = \kop(\al,\be)$. Also, without loss of generality,  we may assume that $\pi(\al) > k + 1$. Then, we distinguish the following cases.

\newcases
\begin{case}
  $m < k$.
\end{case}
As $m+1\leq k$ and $k = \kop(\al,\be)$, we have that $I(\al,m+1) = I(\be,m+1)$. Then as $m = \kop(\al,\zeta)$, we infer that $\zeta\not\in I(\be,m+1)$. Thus, $m = \kop(\be,\zeta)$. We see that since $\pi(s) < \pi(t)$, $\lop(\be) > m+1$.

  Since  $t \prec_p v$ and $p$ satisfies $(P4)$, there is a walk $\langle t, t'_{\lop(\be)},t'_{\lop(\be)-1},\dots,t'_{m+1},$ $v_{m+1},\dots , v_{\lop(\zeta)-1}, v \rangle$ from $t$ to $v$ in $p$. Then, $\langle s, s_{\lop(\al)}, s_{\lop(\al)-1},\dots, s_{k+1},$ $t'_k,t'_{k-1},$ $\dots, t'_{m+1} ,v_{m+1},$ $\dots , v_{\lop(\zeta)-1}, v \rangle$ is a walk from $s$ to $v$ in $q$.

\begin{case}
  $m = k$ and  $\zeta\not\in I(\be,k+1)$.
\end{case}

    We may suppose that $\lop(\be) > k+1$. Otherwise, the considerations are similar.  As $t \prec_p v$, by condition $(P4)$ for $p$, there is a walk $\langle t, t'_{\lop(\be)}, t'_{\lop(\be)-1},\dots,t'_{m+1},$ $v_{m+1},\dots , v_{\lop(\zeta)-1}, v \rangle$ from $t$ to $v$ in $p$. Then, the sequence $\langle s, s_{\lop(\al)}, s_{\lop(\al)-1},\dots,s_{k+1},$ $v_{k+1}, v_{k+2},\dots,$ $v_{\lop(\zeta)-1},v \rangle$ is a walk from $s$ to $v$ in $q$.

      \begin{case}
        $m = k$ and  $\zeta\in I(\be,k+1)\setminus I(\be,\lop(\be))$.
      \end{case}
      
      Let $l = \kop(\be,\zeta)$. Since $\zeta\in I(\be,k+1)$, we see that $k < l$. We may assume that $\lop(\be) > l+1$.  As $t \prec_p v$, there is a walk $\langle t, t'_{\lop(\be)}, t'_{\lop(\be)-1},\dots,t'_{l+1},$ $v_{l+1}, \dots, v_{\lop(\zeta)- 1}, v \rangle$ from $t$ to $v$ in $p$. Then, the sequence $\langle s, s_{\lop(\al)}, s_{\lop(\al)-1},\dots,s_{k+1},$ $t_{k+1}, \dots, $ $ t_l, v_{l+1}, v_{l+2},\dots, v_{\lop(\zeta)-1}, v \rangle$ is a walk from $s$ to $v$ in $q$.

      \begin{case}
        $m = k$ and  $\zeta\in I(\be,\lop(\be))$.
      \end{case}

      Let $l = j(\be,\zeta)$.  As $t \prec_p v$, there is a walk $\langle t, t_{l+1}, v_{l+1},v_{l+2}, \dots , v_{\lop(\zeta)-1},  v \rangle$ from $t$ to $v$ in $p$. Then,  the sequence $\langle s, s_{\lop(\al)}, s_{\lop(\al)-1},\dots,s_{k+1},$ $t_{k+1}, \dots,$ $ t_{\lop(\beta)-1}, t, v_{l+1},\dots, v_{\lop(\zeta)-1},  v \rangle$ is a walk from $s$ to $v$ in $q$.

        \section{\bf Proof of Lemma 2.10.}\label{sc:Lm.2.10}
        
        \vspace{2mm} Actually, we will prove that  ${\mathbb P}_{\delta}$  satisfies  Property K.

      \vspace{2mm} In what follows, we should keep in mind the following definitions. If $v\in {\mathbb B }_{\langle \al,\zeta \rangle}$ for $\al\in \tilde{L}^{\lambda}_{\om_2}$ and $0 < \zeta < \lambda$, we put $o(v) = \al$. 
      
      \vspace{1mm}
      For every $\al < \de$ we write $Y_{\al} = \{y\in Y : \pi(y) = \al \}$. 
      
      \vspace{1mm}
      If $s\in Y$, we write $\lop(s) = \lop(\pi(s))$.

\vspace{1mm}
       In order to prove Lemma 2.10, assume that $R = \{p_{\nu} : \nu < \om_1 \} \subset P_{\de}$ with $p_{\nu}\neq p_{\mu}$ for $\nu < \mu < \om_1$. For each $\nu < \om_1$ we write $p_{\nu} = \langle x_{\nu},\preceq_{\nu},i_{\nu} \rangle$ and $x_{\nu} = \{s_{\nu,i} : i < n_{\nu}\}$. By thinning out $R$ by means of standard combinatorial arguments, we can assume the following:

      \begin{enumerate}[(A)]
      \item $\{x_{\nu} : \nu < \om_1 \}$ forms a $\Delta$-system with kernel $x$ such that:

      \begin{enumerate}[(a)]
      \item $|x_{\nu}| = |x_{\mu}|$ for $\nu < \mu < \om_1$.

      \item For every $\al <\de$, either $x_{\nu}\cap Y_{\al} \cap {\mathbb B}_1 = x \cap Y_{\al} \cap {\mathbb B}_1$ for each $\nu < \om_1$ or there is at most one $\nu < \om_1$ such that $x_{\nu}\cap Y_{\al} \cap {\mathbb B}_1 \neq \emptyset$.

      \end{enumerate}

      \item If $H_{\nu} = \{t\in {\mathbb B} : x_{\nu}\cap {\mathbb B}_t \neq \emptyset \}$ for each $\nu < \om_1$ then the following conditions hold:

    \begin{enumerate}[(a)]
    \item $\{H_{\nu} : \nu < \om_1 \}$ forms a $\Delta$-system with  $|H_{\nu}| = |H_{\mu}|$ for $\nu < \mu < \om_1$.

    \item For every $\al \in \tilde{L}^{\lambda}_{\om_2}$ and every $ 0 < \zeta < \lambda$, if $x\cap {\mathbb B}_{\langle \al, \zeta \rangle} \neq \emptyset$ then $x_{\nu} \cap {\mathbb B}_{\langle \al, \zeta \rangle} = x \cap {\mathbb B}_{\langle \al, \zeta \rangle} $ for each $\nu < \om_1$.
        
    \item  For every $\al \in \tilde{L}^{\lambda}_{\om_2}$, every $ 0 < \zeta < \lambda$ and every $\nu < \om_1$, if $x_{\nu}\cap  {\mathbb B }_{\langle \al,\zeta \rangle}\neq \emptyset$ and $x_{\nu}\cap  {\mathbb B}_1 \cap Y_{\ga_{\al} + 1} = \emptyset$ then $x_{\mu}\cap  {\mathbb B}_1 \cap Y_{\ga_{\al} + 1} = \emptyset$ for every $\mu < \om_1$.

      \end{enumerate}

\item For each ${\nu}<{\mu}<\omega_1$ there is an isomorphism
$g=g_{{\nu},{\mu}}:\<x_{\nu},\preceq_{\nu}, i_{\nu}\>\to
\<x_{\mu},\preceq_{\mu},i_{\mu}\>$ such that  for every $s,t\in x_{\nu}$ the following conditions hold:

\begin{enumerate}[(1)]

\item $g\restriction x =\operatorname{id}_x$,
\item  $g(s_{\nu,i}) = g(s_{\mu,i})$ for all $i < n_{\nu} = n_{\mu}$,
\item  $\pi_{\mathbb B}(s) = \pi_{\mathbb B}(t)$ iff $\pi_{\mathbb B}(g(s)) = \pi_{\mathbb B}(g(t))$,
\item  $\pi_{\mathbb B}(s) =  1$ iff $\pi_{\mathbb B}(g(s)) = 1$,
\item $\pi(s) = \pi_{-}(s)$ iff $\pi(g(s)) = \pi_{-}(g(s))$,
\item  $\pi_{-}(s) = \pi_{-}(t)$ iff $\pi_{-}(g(s)) = \pi_{-}(g(t))$,
\item $i_{\nu}\{s,t\} = i_{\mu}\{s,t\}$ for all $\{s,t\}\in [x]^2$,
\item $\lop(s) = \lop(g(s))$.

 \end{enumerate}

     \end{enumerate}

     \vspace{1mm}
     Note that in order to obtain $(C)(7)$ we use condition $(P3)(d)$ and the fact that $h\{s,t\}$ is countable for all $\{s,t\}\in [Y]^2$.

     \vspace{1mm}
     Let $n = n_{\nu}$ for $\nu < \om_1$ and let $M = \{i\in n : s_{\nu,i}\not\in U_2 \}$. Then, we may also assume the following:

      \begin{enumerate}[(D)]
      \item There is a partition $M = K \cup^* D \cup^* L \cup^* F$ such that for each $\nu < \mu < \om_1$ the following holds:

      \begin{enumerate}[(1)]
      \item $\forall i\in K$ $s_{\nu,i}\in x$, and so $s_{\nu,i} = s_{\mu,i}.$

      \item $\forall i\in D$ $\pi_B(s_{\nu,i}) = 1$ and $\pi(s_{\nu,i}) \neq \pi(s_{\mu,i}).$

      \item $\forall i\in L$  $\pi_B(s_{\nu,i}) \neq 1$, $\pi_B(s_{\nu,i}) \neq \pi_B(s_{\mu,i})$ but $\pi(s_{\nu,i}) = \pi(s_{\mu,i}).$

       \item $\forall i\in F$ $\pi_B(s_{\nu,i}) \neq 1$  and $\pi(s_{\nu,i}) \neq \pi(s_{\mu,i}).$

      \end{enumerate}

       \end{enumerate}

     \vspace{1mm} The proofs of the next two claims are similar to the proofs of  \cite[Lemmas 4 and 6]{ma1}. We show these claims for completeness.

     \begin{claim} Assume that $\nu < \mu  < \om_1$ and  $u,s\in x_{\nu}\cap ({\mathbb B}_1\cup U_1)$ with  $u \prec_{\nu} s$. Let $u'  = g_{\nu,\mu}(u)$  and $s'  = g_{\nu,\mu}(s)$. 
     
     \vspace{1mm} (a) If   $I\in {\mathbb I}$ with  $\pi(u),\pi(u'),\pi(s)\in I$ then $\pi(s')\in I$.

     \vspace{1mm} (b) If   $I\in {\mathbb I}$ with  $\pi(u),\pi(s),\pi(s')\in I$ then  $\pi(u')\in I$.
     \end{claim}

 \vspace{2mm}
 \noindent
{\bf Proof.}  We show condition $(a)$. The proof of $(b)$ is similar. So, assume that  $I\in {\mathbb I}$ with  $\pi(u),\pi(u'),\pi(s)\in I$. Suppose that $I\in \ical_n$. And assume on the contrary that 
      $\pi(s')\not\in I$. Let $J= I(\pi(s'),n)$.  Since $\pi(s')\not\in I$, we have $I\neq J$. Since $u \prec_{\nu} s$, we see that $u' \prec_{\mu} s'$. Note that $J^-$ is in the sequence $w(\pi(u'),\pi(s'))$. Then as $u'\prec_{\mu} s'$, by condition $(P4)$, there is $v'\in x_{\mu}$ with $\pi(v') = J^-$ such that $u' \prec_{\mu} v' \preceq_{\mu} s'$. As $\pi(v') = J^-$, we infer that  $\lop(v') \leq n$. Let $v = g^{-1}_{\nu,\mu}(v')$. Then $u \prec_{\nu} v \preceq_{\nu} s$ and, by condition $(C)(8)$, $\lop(v) = \lop(v') \leq n$.  But this is impossible, because $\pi(u),\pi(v),\pi(s)\in I \in \ical_n$.    $\square$
      
      \vspace{1mm} Note that if in the statement of Claim 3.1, we have that $u\prec_{\nu}s$ and  $u\in x$, then  $u' = u$, and so if $I\in {\mathbb I}$ with  $\pi(u),\pi(s)\in I$, we deduce that $\pi(s')\in I$. This particular case of Claim 3.1 will be frequently used.

      \vspace{2mm} We write $c_{\nu} = \pi[(x_{\nu}\setminus x)\cap ({\mathbb B}_1\cup U_1)]$ for $\nu < \om_1$. Note that $|c_{\nu}| = |c_{\mu}|$
      for $\nu < \mu < \om_1$. Suppose that $A$ is an uncountable subset of $\om_1$. Then, we say that $\{p_{\nu} : \nu\in A \}$ is  {\em admissible}, if for every $I\in \ical_n$ such that $\cf (I^+) = \om_2$ and $\{c_{\nu}\cap E(I) : \nu \in A \}$ is uncountable, we have that $\{c_{\nu}\cap E(I) : \nu \in A \}$ is adequate for the function $G_I$.

      \begin{claim} There is an uncountable subset $A$ of $\om_1$ such that $\{p_{\nu} : \nu\in A \}$ is admissible.
     \end{claim}

 \vspace{2mm}
 \noindent
{\bf Proof.} Let $\{i_1,\dots ,
       i_r \} = M \cap (D\cup F)$ where $i_1 < i_2 < \dots < i_r$.
      Then, we define a sequence $A_0 \supset A_1 \supset \dots \supset A_r$ of uncountable subsets of $\om_1$ as follows. We put $A_0 = \om_1$. Assume that $1\leq k \leq r$. We choose an $n < \om$ and an interval $J_k\in \ical_n$ such that $\{\pi(s_{\nu,i_k})\in J_k : \nu\in A_{k-1} \}$ is uncountable but for every $J\in {\mathbb I}$ with $J\subsetneqq J_k$, $\{\pi(s_{\nu,i_k})\in J : \nu\in A_{k-1} \}$ is countable.
       Note that this interval exists, because otherwise we would find an infinite strictly decreasing sequence of ordinals in $\de$. Now, we put $A'_k = \{\nu\in A_{k-1}: \pi(s_{\nu,i_k}) \in J_k \}$. Then, we define $A_k$ as follows. If $|E(J_k)| = \om_2$ and $\{ \pi(s_{\nu,i_k}) : \nu \in A'_k \}\subset E(J_k)$, then $A_k$ is an uncountable subset of $A'_k$ such that $\{c_{\nu} \cap E(J_k) : \nu \in A_k \}$ is adequate for the function $G_{J_k}$. Otherwise, $A_k = A_{k-1}$.

      \vspace{1mm} We define $A = A_r$. To verify that $A$ is as required, consider $I\in \ical_n$ such that $\cf (I^+) = \om_2$ and $\{c_{\nu} \cap E(I) : \nu \in A \}$ is uncountable. Then, for some $1\leq k \leq r$, $\{\pi(s_{\nu,i_k}) : \nu\in A\}\cap E(I)$ is uncountable. Hence, by the construction of $A$, we have that $I = J_k$ and  $\{c_{\nu} \cap E(I) : \nu \in A \}$ is adequate for $G_I$.    $\square$

       \vspace{2mm} So, assume that $A$ is an uncountable subset of $\om_1$ such that $\{p_{\nu} : \nu \in A \}$ is admissible. Let $\nu,\mu \in A$ with $\nu < \mu$.
Our aim is to show that $p_{\nu}$ and $p_{\mu}$ are compatible in ${\mathbb P}_{\de}$. For this, we define $p = \langle Z,\preceq, i \rangle$ as follows. We put $Z = x_{\nu} \cup x_{\mu}$. If $s,t\in Z$, we put $s\preceq t$ iff $s \preceq_{\nu} t$ or $s \preceq_{\mu} t$ or  there is $u\in x$ such that either $s\preceq_{\nu} u \preceq_{\mu} t$ or $s\preceq_{\mu} u \preceq_{\nu} t$. And for $\{s,t\}\in [Z]^2$ we define $i\{s,t\}$ as follows: $i\{s,t\} = i_{\nu}\{s,t\}$ if $s,t\in x_{\nu}$;  $i\{s,t\} = i_{\mu}\{s,t\}$ if $s,t\in x_{\mu}$;  and if  $s\in x_{\nu}\setminus x_{\mu}$ and $t\in x_{\mu}\setminus x_{\nu}$, we define $i\{s,t\} = \{u\in Z : u\preceq s,t \mbox{ and } \pi(u)\in h\{s,t\} \}$ if $s,t$ are not comparable in $p$, $i\{s,t\} = \{s\}$ if $s\prec t$, and $i\{s,t\} = \{t\}$ if $t\prec s$. Note that, by condition $(C)(7)$, the definition of $i$ is meaningful.

  \vspace{1mm} We have to show that $p\in P_{\de}$. It is easy to check that $p$ satisfies conditions $(P1)$, $(P2)$ and $(P3)(a)-(e)$. So, we prove that $p$ satisfies condition $(P3)(f)$. To verify this point, consider $s,t\in Z$ with $s\neq t$. First, suppose that $s,t\in x_{\nu}$. Let $u\in Z$ be such that $u\preceq s,t$. If $u\in x_{\nu}$, we are done. So, suppose that $u\in x_{\mu}\setminus x$.  As $u\prec s,t$, there are $w_1,w_2\in x$ such that $u\prec_{\mu} w_1 \preceq_{\nu} s$ and  $u\prec_{\mu} w_2 \preceq_{\nu} t$. We may assume that $w_1\neq w_2$. Let $w\in i_{\mu}\{w_1,w_2\}$ be such that $u\prec_{\mu} w$. Note that $w\in x$ by condition  $(C)(7)$.  Hence, there is $v\in i_{\nu}\{s,t\}$ such that $w\preceq_{\nu} v$, and so $u\prec v$. If $s,t\in  x_{\mu}$,  we can use a parallel argument.

  \vspace{1mm} Now, assume that $s\in x_{\nu}\setminus x$ and $t\in x_{\mu}\setminus x$ are compatible but not comparable in $p$. Suppose that $\pi_{-}(s) \leq \pi_{-}(t)$. Otherwise, the argument is parallel. Let $u\in Z$ be such that $u \prec s,t$.   We put $a = \pi[x]$, $a_{\nu} = \pi[x_{\nu}]$ and  $a_{\mu} = \pi[x_{\mu}]$. We distinguish the following cases.

  \vspace{1mm} \noindent {\bf Case 1}. $\pi_{-}(s) = \pi_{-}(t)$.

  \vspace{1mm} Note that by condition   $(A)(b)$, either $s\in U$ or $t\in U$. So, 
we may assume that $s\in  U$ and  $t\in {\mathbb B}_1 \cup U$. Note that by condition $(B)(a)-(b)$, we have   $\pi_{\mathbb B}(s) \neq \pi_{\mathbb B}(t)$. And by condition $(P2)(b)$, $u\in {\mathbb B}_1$. Let $I = J(\pi_{-}(s))$. We distinguish the following subcases.

  \vspace{1mm} \noindent {\bf Case 1.1}. $\pi(u)\in a\setminus I$.

   Since $u \prec_{\nu} s$, by condition $(P4)$, there is $v_1\in x_{\nu}\cap {\mathbb B}_1$ with $\pi(v_1) = I^-$ such that $u \prec_{\nu} v_1 \prec_{\nu} s$. Analogously,  there is $v_2\in x_{\mu}\cap {\mathbb B}_1$ with  $\pi(v_2) = I^-$ such that $u \prec_{\mu} v_2 \prec_{\mu} t$. Now, by conditions $(P3)(b)$ and  $(A)(b)$, we deduce that $v_1 = v_2 \in x$. Thus, since $I^-\in E(J(\pi_{-}(s))) = h\{s,t\}$, we are done.

   \vspace{1mm} \noindent {\bf Case 1.2}. $\pi(u)\in a\cap I$.

   \vspace{1mm} Let $n = \kop(\pi(u),\pi_{-}(s))$. Note that, by condition $(*)(2)$, we have that $\lop(\pi_{-}(s)) = \lop(o(\pi_{-}(s))) + 1$, because if we put $\beta = o(\pi_{-}(s))$ then $\gamma_{\beta} = \pi_{-}(s)$.  Then since $\pi(u)\in  I = J(\pi_{-}(s))$, we have that  $n = \lop(o(\pi_{-}(s)))$, $J(o(\pi_{-}(s)))$ is the interval where $\pi(u)$ and $\pi_{-}(s)$ separate, $I(\pi(u),n+1) = I$ and $I(\pi_{-}(s),n+1) = [\pi_{-}(s),\pi_{-}(s) +2)$ by Definition 2.1.

   \vspace{1mm} Let $E(I) = \langle \epsilon_k : k < \om \rangle$. We  assume that $\pi(u)\not\in E(I)$.  Otherwise, the case is obvious. Let $m$ be the least natural number $n$ such that $\pi(u) < \epsilon_n$. Note that, by Proposition \ref{Pr.2.6},  $\epsilon_m$ is in the sequence $w(\pi(u),\pi(s))$ and also in the sequence $w(\pi(u),\pi(t))$. Since $u \prec_{\nu} s$ and $u \prec_{\mu} t$, by condition $P(4)$, there are $v_1\in x_{\nu}\cap {\mathbb B}_1$ and $v_2\in x_{\mu}\cap {\mathbb B}_1$ with $\pi(v_1) = \pi(v_2) = \epsilon_m$ such that $u \prec_{\nu} v_1 \prec_{\nu} s$ and $u \prec_{\mu} v_2 \prec_{\mu} t$. So, we have $v_1 = v_2\in x$ by conditions $(A)(b)$ and $(P3)(b)$. Hence, as $\pi(v_1)\in E(J(\pi_{-}(s))) = h\{s,t\}$, we are done.

   \vspace{1mm} \noindent {\bf Case 1.3}. $\pi(u)\in (a_{\nu}\setminus a)\cap I$.

   \vspace{1mm} Since $u \prec t$, there is $w\in x$ such that $u \prec_{\nu} w \prec_{\mu} t$.  Let $E(I) = \langle \epsilon_k : k < \om \rangle$. We may assume that $\pi(u)\not\in E(I)$. Let $m$ be the least natural number $n$ such that $\pi(u) < \epsilon_n$. First, assume that $\epsilon_m \leq \pi(w)$. By  condition $(P4)$, there are $v_1,v_2\in x_{\nu} \cap {\mathbb B}_1$ with $\pi(v_1) = \pi(v_2) = \epsilon_m$ such that  $u \prec_{\nu} v_1 \prec_{\nu} s$ and $u \prec_{\nu} v_2 \preceq_{\nu} w$. By $(P3)(b)$, $v_1 = v_2$. So as $v_1 \prec s,t$ and $\pi(v_1) \in E(J(\pi_{-}(s))) = h\{s,t\}$, we are done.

    \vspace{1mm} Now, assume that $\epsilon_m > \pi(w)$.  By  condition $(P4)$, there are $v_1\in x_{\nu}\cap {\mathbb B}_1$ and $v_2\in x_{\mu}\cap {\mathbb B}_1$
    with $\pi(v_1) = \pi(v_2) = \epsilon_m$ such that  $u \prec_{\nu} v_1 \prec_{\nu} s$ and $w \prec_{\mu} v_2 \prec_{\mu} t
    $. Therefore, $v_1 = v_2 \in x$, and so $v_1$ is as required.

     \vspace{2mm} Also, if either $\pi(u)\in a_{\nu}\setminus (a\cup I)$ or $\pi(u)\in a_{\mu}\setminus (a\cup I)$, we can proceed by means of an argument similar to the one given in Case 1.1. And if $\pi(u)\in (a_{\mu}\setminus a)\cap I$, we can proceed as in Case 1.3.

     \vspace{2mm}  \noindent {\bf Case 2}. $\pi_{-}(s) < \pi_{-}(t)$.

     \vspace{1mm} Let $I \in \ical_k$ be the interval where $\pi_{-}(s)$ and $\pi_{-}(t)$ separate. Let $J = I(\pi_{-}(s),k+1)$ and $K = I(\pi_{-}(t),k+1)$.

     \vspace{1mm} Note that, by condition $(B)(a)-(b)$, if either $s\not\in {\mathbb B}_1$ or $t\not\in {\mathbb B}_1$, we see that $\pi_{{\mathbb B}}(s) \neq \pi_{{\mathbb B}}(t)$. So, by condition 
     $(P2)(b)$, $u\in {\mathbb B}_1$.

     \vspace{1mm} We assume that $u\in x_{\nu}$. Otherwise, the argument is similar.  If $u\not\in x$, we pick an element $w\in x$ such that $u \prec_{\nu} w \prec_{\mu} t$. Note that $w\in {\mathbb B}_1$ by condition $(B)(a)-(b)$.

     \vspace{1mm} Now, we define the elements $s^*$ and $t^*$ as follows.
      If $s\in U_2$, by using condition  $(P2)(c)$, we pick $ s^* \in x_{\nu}\cap U_1$ with $\pi_{{\mathbb B}}(s^*) = \pi_{{\mathbb B}}(s)$ such that $u \prec_{\nu} s^* \prec_{\nu} s$,
       and otherwise we define $s^* = s$. Now, assume that $t\in U_2$. If $u\in x$, we pick  $t^* \in x_{\mu}\cap U_1$ with  $\pi_{{\mathbb B}}(t^*) = \pi_{{\mathbb B}}(t)$ such that
       $u \prec_{\mu} t^* \prec_{\mu} t$ by using again $(P2)(c)$. And if $u\not\in x$, we pick
      $t^*\in x_{\mu}\cap U_1$ with $\pi_{\mathbb B}(t^*) = \pi_{\mathbb B}(t)$ such that $w\prec_{\mu} t^* \prec_{\mu} t$. If $t\not\in U_2$,  we put $t^* = t$.
      Note that, by condition $(B)(b)$,  $s^*\in x_{\nu}\setminus x$ and $t^*\in x_{\mu}\setminus x$. And clearly, $\pi(s^*)\in J$ and $\pi(t^*)\in K$. We have that $s^*\not\prec t^*$. To check this point, assume on the contrary that  $s^* \prec t^*$. Note that since $\pi(s^*) < \pi(t^*)$, by condition $(P2)(b)$, it follows that $s\in {\mathbb B}_1$, and so $s^* = s$.  Let $v\in x$ be such that $s^* \prec_{\nu} v \prec_{\mu} t^*$. So, we see that $s^* = s \prec_{\nu} v \prec_{\mu} t^* \preceq_{\mu} t$, which contradicts the assumption that $s$ and $t$ are incomparable in $p$.

       \vspace{1mm}  Let $\tilde{s} = g_{\nu,\mu}(s^*)$ and $\tilde{t} = g_{\nu,\mu}^{-1}(t^*)$. Recall that if  $v\in {\mathbb B }_{\langle \al,\zeta \rangle}$ for $\al\in \tilde{L}^{\lambda}_{\om_2}$ and $0 < \zeta < \lambda$, we have defined $o(v) = \al$. 

      \begin{claim} Assume that $\pi(u)\in I$ with $\pi(u) > I^-$.

       \vspace{1mm}
      (a) If $s\in U$ then $I\supsetneqq J(o(s))$.

       \vspace{1mm}
      (b) If $t\in U$ then $I\supsetneqq J(o(t))$.

      \end{claim}
      
       \noindent {\bf Proof}. First, we prove $(a)$. So, suppose that $s\in U$. Let $n\in \om$ be such that $J(o(s)) = J(o(s^*)) \in \ical_n$. So, $I(\pi(s^*),n) =   J(o(s^*))$ and $I(\pi(s^*),n +1) = [\pi(s^*),\pi(s^*) + 2)$ by condition $(*)(4)$,  because if we put $\beta = o(s^*)$ then $\gamma_{\beta} = \pi(s^*)$.  Since $u\in I$, we see that  $I\supsetneqq I(\pi(s^*), n+1)$, and so $I\supseteq I(\pi(s^*),n) = J(o(s^*))$. We have to show that $I\supsetneqq J(o(s))$. So, assume on the contrary that $I = J(o(s))$. Then, $J = [\pi(s^*),\pi(s^*) + 2)$. Since $\pi(s^*) < \pi(t^*)$, we infer that  $t^*\in {\mathbb B}_1$ by $(B)(a)-(b)$ and $(*)(3)$, and so  $\tilde{t}\in  {\mathbb B}_1$ by condition $(C)(4)$. Also,  $\tilde{s}\in U$ again by condition $(C)(4)$. By using Claim 3.1, we can check that $\pi(\tilde{s}), \pi(\tilde{t})\in I$. For this, if $u\in x$, use the fact that $u\prec_{\nu} s^*$ and  $u\prec_{\mu} t^*$. Now, suppose that $u\not\in x$. Since $w \prec_{\mu} t^*$ and $\pi(t^*)\in I\setminus J$, we see that  $\pi(\tilde{t})\in I\setminus J$. Let $\tilde{u} = g_{\nu,\mu}(u)$. As $u\prec_{\nu} w$, we see that  $\tilde{u} \prec_{\mu} w$, and so we infer that $\pi(\tilde{u})\in I$. Hence, as $u \prec_{\nu} s^*$ and $\pi(u),\pi(\tilde{u}), \pi(s^*)\in I$, we deduce that $\pi(\tilde{s})\in I$. Therefore,  $\pi(\tilde{s}) < \pi(s^*)$ by  $(B)(a)-(b)$ and  $(*)(3)$. Since $\tilde{t}\in {\mathbb B}_1$, by $(P2)(b)$, it is impossible that $s^*  \prec_{\nu}  \tilde{t}$. And if $\tilde{t} \prec_{\nu}  s^*$ then $t^* \prec_{\mu} \tilde{s}$, and so $\pi(t^*) < \pi(\tilde{s})$, which is impossible because $\pi(\tilde{s}) < \pi(s^*)$ and $\pi(s^*) < \pi(t^*)$. Hence, $s^*$ and $\tilde{t}$ are incomparable in $p_{\nu}$. Note that since $\pi(\tilde{s}) < \pi(s^*) < \pi(t^*)$, by condition $(C)(6)$, we deduce that $\pi(\tilde{t}) \neq \pi(s^*)$, and thus $h\{s^*,\tilde{t}\} =  G\{\pi(s^*),\pi(\tilde{t})\}$.
       Since $u \prec_{\nu} s^*, \tilde{t}$, there is $v\in i_{\nu}\{s^*,\tilde{t}\}$ such that $u\preceq_{\nu} v$, and so we infer from $(P3)(d)$ that $\pi(v)\in h\{s^*,\tilde{t}\} = G\{\pi(s^*),\pi(\tilde{t})\} = G_I\{\pi(s^*),\pi(\tilde{t})\} \cup \{I^-\} = \{I^-\}$, because $I^-$ is the first element of $E(I)$ and $\pi(s^*)$ is its  second element by Definition 2.1.  Hence, $\pi(u) \leq I^{-}$, which is impossible. Thus, we have shown that $I\supsetneqq J(o(s))$. 

        \vspace{1mm} Next, we show $(b)$. In this case, the argument is simpler.  Assume that $t\in U$. Thus,  $\tilde{t}\in U$.  Let $n\in \om$ be such that $J(o(t)) = J(o(t^*)) \in \ical_n$. So, $I(\pi(t^*),n) =   J(o(t^*))$ and $I(\pi(t^*),n +1) = [\pi(t^*),\pi(t^*) + 2)$ by $(*)(4)$. Since $\pi(s^*) < \pi(t^*)$, we see that $I\supsetneqq I(\pi(t^*), n + 1)$,  and hence $I\supset J(o(t^*))$.  Suppose that $I = J(o(t^*))$. Then,  $J = J(\pi(t^*))$ and $K = [\pi(t^*),\pi(t^*) + 2)$. Since $u \prec t^*$ and $\pi(t^*)\not\in J$, we infer from Claim 3.1 that $\pi(\tilde{t})\in I\setminus J$, and so $\pi(t^*) \leq \pi(\tilde{t}) < o(t^*)$, which is impossible by $(B)(a)-(b)$ and condition $(*)(3)$.  $\square$

        \vspace{1mm} Now, we distinguish the following subcases.

     \vspace{1mm} \noindent {\bf Case 2.1}. $\pi(u)\in a\setminus  I$.

   \vspace{1mm} We prove that $I^- < \pi(s^*)$. First, suppose that $s\in U$. Note that  $I\supsetneqq [\pi(s^*),\pi(s^*) + 2)$ by condition $(B)(a)-(c)$, and so $I^- < \pi(s^*)$  by $(*)(8)$. Now, suppose that $s\in {\mathbb B}_1$ and assume on the contrary that  $I^- = \pi(s)$. Then as $u \prec_{\mu} t$, by condition $(P4)$, there is $v\in x_{\mu} \cap {\mathbb B}_1$ with $\pi(v) = I^-$ such that  $u \prec_{\mu} v \prec_{\mu} t$. Hence, $s = v \in x$ by $(P3)(b)$, which contradicts the fact that $s\in x_{\nu}\setminus x$.

So, $I^- < \pi(s^*)$. Since $u \prec_{\nu} s$, by condition $(P4)$, there is $v_1\in x_{\nu}\cap {\mathbb B}_1$ with $\pi(v_1) = I^-$ such that $u \prec_{\nu} v_1 \preceq_{\nu} s$. Analogously, as $u \prec_{\mu} t$, there is $v_2\in x_{\mu}\cap {\mathbb B}_1$ with $\pi(v_2) = I^-$ such that $u \prec_{\mu} v_2 \prec_{\mu} t$. Thus, $v_1 = v_2$ and $v_1 \in x$. So, as $I^- \in G\{\pi_{-}(s),\pi_{-}(t)\} = h\{s,t\}$, we are done.

     \vspace{2mm} \noindent {\bf Case 2.2}. $\pi(u)\in a \cap I$.

     \vspace{1mm} We assume that $I^- < \pi(u)$. Otherwise, the case is obvious.  We consider the following subcases.

     \vspace{2mm} \noindent {\bf Case 2.2.1}. $\pi(u)\in J$.

      \vspace{1mm} We show that this case is impossible. As $u\prec_{\nu} s^*$ and $u\prec_{\mu} t^*$,  we deduce from Claim 3.1 that $\pi(\tilde{s})\in J$ and $\pi(\tilde{t})\in I\setminus J$.

      First, assume that $s^* \prec_{\nu} \tilde{t}$. Then, as $s^* \prec_{\nu} \tilde{t}$ and $\tilde{s} \prec_{\mu} t^*$, by $(P4)$, $(P3)(b)$ and $(A)(b)$, we deduce that there is $v\in x\cap {\mathbb B}_1$ with $\pi(v) = J^+$ such that $s^* \prec_{\nu} v \preceq_{\nu} \tilde{t}$ and $\tilde{s}\prec_{\mu} v \preceq_{\mu} t^*$, and thus $s^* \prec t^*$, which is impossible.

      Now, suppose that $s^*\not\prec_{\nu} \tilde{t}$. Since $\pi(\tilde{t})\not\in J$, we see that $s^*$ and $\tilde{t}$ are incomparable in $p_{\nu}$. Let $K_0 = I(\pi(\tilde{t}), k +1)$. Note that as $u\prec_{\mu} t^*$, we have $u\prec_{\nu} \tilde{t}$.
        Then since $u \prec_{\nu} s^*, \tilde{t}$, there is $v\in i_{\nu}\{s^*,\tilde{t}\}$ such that $u\preceq_{\nu} v$,  and so we infer from $(P3)(d)$ that $\pi(v)\in h\{s^*, \tilde{t}\} = G_I\{J^-,K^-_0\} \cup \{I^-\}$. Note that since  $G_I \{J^-,K^-_0\} \subset J^-$, $\pi(u)\in J$ and $\pi(u) > I^-$, we infer that $\pi(v) < \pi(u)$, which contradicts the fact that $u\preceq_{\nu} v$.

     \vspace{2mm} \noindent {\bf Case 2.2.2}. $\pi(u)\not\in J$.

      \vspace{1mm} Since $u\prec_{\nu} s^*$ and $u\prec_{\mu} t^*$, by $(P4)$, there are $v_1\in x_{\nu}\cap Y_{J^-}$ and $v_2\in x_{\mu}\cap Y_{K^-}$ such that $u \prec_{\nu} v_1 \preceq_{\nu} s^*$ and  $u \prec_{\mu} v_2 \preceq_{\mu} t^*$. Clearly, if $s\in {\mathbb B}_1$ then $v_1\in {\mathbb B}_1$. And
       note that if $s\in U$  we have that $I\supsetneqq J(o(s))$ by Claim 3.3(a), and hence $J \supsetneqq [\pi(s^*), \pi(s^*) + 2)$ by $(*)(7)$, and so $\pi(s^*) > J^-$ by $(*)(8)$. Analogously, if $t\in {\mathbb B}_1$ then $v_2\in {\mathbb B}_1$, and if $t\in U$  we have that $I\supsetneqq J(o(t))$ by Claim 3.3(b), and hence $K \supsetneqq [\pi(t^*), \pi(t^*) + 2)$ by $(*)(7)$, and so $\pi(t^*) > K^-$ by $(*)(8)$. Therefore, $v_1,v_2\in {\mathbb B}_1$. Recall that our aim is to show that there is $v\in Z$ such that $u\preceq v \prec s,t$ and $\pi(v)\in h\{s,t\}$. We distinguish the following subcases.

\vspace{2mm} \noindent {\bf Case 2.2.2.1}. $v_1\in x$.

      \vspace{1mm} Since  $u \prec_{\nu} v_1 \prec_{\nu} s^*$ and $v_1\in J$, by Case 2.2.1,  it is impossible that  $v_1 \prec_{\mu} t^*$. Hence,
     $v_1$ and $t^*$ are incomparable in $p_{\mu}$.
       We have that $h\{v_1,t^*\} = G_I\{J^-,K^-\} \cup \{I^-\}$.
 Then as $u \prec_{\mu} v_1,t^*$, there is $v\in i_{\mu}\{v_1, t^*\}$ such that $u \preceq_{\mu} v$. Therefore, we have $v\prec_{\mu} v_1 \prec_{\nu} s^*$ and
$v\prec_{\mu} t^*$, and thus $u\preceq v \prec s,t$. By $(P3)(d)$, we infer that $\pi(v)\in h\{v_1,t^*\}$. Note that as $v_1$ and $t^*$ are incomparable in $p_{\mu}$,
we have $v\neq v_1$. Then as $I(\pi(s^*),k+1) = I(\pi(v_1),k+1)$, we deduce that $h\{v_1,t^*\} = h\{s^*,t^*\}$, and hence $\pi(v)\in h\{s^*,t^*\} = h\{s,t\}$.

\vspace{2mm} \noindent {\bf Case 2.2.2.2}. $v_2\in x$.

      \vspace{1mm} As $s^* \not\prec t^*$, we see that $s^*$ and $v_2$ are incomparable in $p_{\nu}$, and then since $u \prec_{\nu} s^*,v_2$ and $I(\pi(v_2),k+1) = I(\pi(t^*),k+1)$, we can proceed as in the previous subcase.

\vspace{2mm} \noindent {\bf Case 2.2.2.3}. $v_1\not\in x$ and $v_2\not\in x$.

 Let $J_0 = I(\pi(u), k+1)$. First, we show that $J^+_0 < J^-$. For this, assume on the contrary that $J^+_0 = J^-$. Since $u \prec_{\mu} v_2$, by using condition $(P4)$,  we infer that there is $v'_1\in x_{\mu}\cap Y_{J^-}$ such that $u \prec_{\mu} v'_1 \prec_{\mu} v_2$. So, by conditions $(A)(b)$ and $(P3)(b)$, we deduce that $v_1 = v'_1 \in x$, which is impossible. Hence,  $J^+_0 < J^-$. Now, since $u \prec_{\nu} v_1$ and $u \prec_{\mu} v_2$, we infer from conditions $(P4)$, $(P3)(b)$ and $(A)(b)$ that there is $v\in x \cap {\mathbb B}_1$ with
 $\pi(v) = J^+_0$ such that $u \prec_{\nu} v$, $v\prec_{\nu} v_1$ and $v\prec_{\mu} v_2$. Note that $\pi(v)\in E(I)$, because $\pi(v) = J^+_0$.  Our purpose is to show that $\pi(v)\in h\{s,t\}$.
For this, we deduce from Claim 3.1 and condition $(C8)$ that $\pi(g_{\nu,\nu'}(v_1))\in E(I)$ for every $\nu'\in A$ with $\nu < \nu'$.
And since $v_1\in {\mathbb B}_1\setminus x$, we infer that for $\nu'\in A$ with $\nu < \nu'$ we have that $g_{\nu,\nu'}(v_1)\in {\mathbb B}_1\setminus x$ by conditions $(C1)$ and $(C4)$,
and hence $\{\pi(g_{\nu,\nu'}(v_1)) : \nu < \nu', \nu'\in A \}$ is uncountable. Recall that we defined $c_{\nu} = \pi[(x_{\nu}\setminus x)\cap ({\mathbb B}_1\cup U_1)]$ for $\nu < \om_1$.
Thus, $\{c_{\gamma}\cap E(I) : \gamma \in A \}$ is uncountable.
 We may assume that $|E(I)| = \om_2$. Otherwise, the argument is simpler. Then as $\{p_{\gamma} : \gamma \in A \}$ is admissible, $\{c_{\gamma}\cap E(I) : \gamma \in A \}$ is adequate for the function $G_I$. We have $\pi(v),\pi(v_1),\pi(v_2)\in E(I)$ and $\pi(v) < \pi(v_1) < \pi(v_2)$. Since $v_1,v_2\not\in x$, we see that $J^-\in a_{\nu}\setminus a$ and $K^-\in a_{\mu}\setminus a$. Then as $\pi(v)\in a$, by adequateness,  we deduce that $\pi(v)\in G_I\{J^-,K^-\} \subset G\{J^-, K^-\} = h\{s^*,t^*\} = h\{s,t\}$.

\vspace{2mm}
      \noindent {\bf Case 2.3}. $\pi(u)\in (a_{\nu}\setminus a)\cap I$.

     \vspace{1mm} We can assume that $\pi(u) > I^-$.  Recall that we have $u \prec_{\nu} s^*$  and $u \prec_{\nu} w  \prec_{\mu}  t^*$ where $w\in x$. Clearly, $\pi(w)\in I$, because $\pi(u),\pi(t^*)\in I$ and $\pi(u) < \pi(w) < \pi(t^*)$. Note that it is not possible that $s^* \prec_{\nu} w$, because $s^*\not\prec t^*$. And if $w\prec_{\nu} s^*$, we can apply Case 2.2 to $s$, $t$ and $w$.  So, we may assume that $s^*$ and $w$ are incomparable in $p_{\nu}$. Assume that $\pi(w)\in J$. Since $w\prec_{\mu} t^*$ and $\pi(t^*)\not\in J$, we infer from Claim 3.1 that $\pi(\tilde{t})\not\in J$. Then as $w\prec_{\nu}\tilde{t}$ and $w\prec_{\mu} t^*$, we deduce from $(P4)$, $(P3)(b)$ and $(A)(b)$ that there is $w'\in x$ with $\pi(w') = J^+$ such that $w\prec_{\nu} w'$,  $w'\prec_{\nu} \tilde{t}$ and $w'\prec_{\mu} t^*$. So,  we can replace $w$ with $w'$, and thus we may assume that $\pi(w)\not\in J$. Let $K_0 = I(\pi(w),k+1)$.

     We prove that $\pi(u)\not\in J$. For this, assume on the contrary that $\pi(u)\in J$.  Then, since $u \prec_{\nu} s^*, w$, there is $v'\in i_{\nu}\{s^*,w\}$ such that $u\preceq_{\nu} v'$,  and so we infer from $(P3)(d)$ that $\pi(v')\in h\{s^*, w\} = G_I\{J^-,K^-_0\} \cup \{I^-\}$. Note that since  $G_I \{J^-,K^-_0\} \subset J^-$, $\pi(u)\in J$ and $\pi(u) > I^-$, we infer that $\pi(v') < \pi(u)$, which contradicts the fact that $u\preceq_{\nu} v'$. Thus, we have that $\pi(u)\not\in J$.

      Next, suppose that $\pi(w)\in K$.   As
      $u \prec_{\nu} s^*,w$, we infer from $(P3)(f)$ and $(P3)(d)$  that there exists $v\in i_{\nu}\{s^*,w\}$ such that $u\preceq v$ and
      $\pi(v)\in h\{s^*,w\} = h\{s^*,t^*\} = h\{s,t\}$. So, as $v \prec s,t$, we are done.

       Now, assume that $\pi(w)\not\in K$. Thus $\pi(t^*)\not\in K_0$ and, by Claim 3.1,  $\pi(\tilde{t})\not\in K_0$. Then as $w\prec_{\nu}  \tilde{t}$ and $w\prec_{\mu} t^*$, by the argument given above, we may assume that  $\pi(w)\in E(I)$.
       By using again condition $(P4)$, it follows that there are $v_1\in x_{\nu}$ with $\pi(v_1) = J^-$ such that $u\prec_{\nu} v_1 \preceq_{\nu} s^*$ and $v_2\in x_{\mu}$ with $\pi(v_2) = K^-$ such that $w\prec_{\mu} v_2 \preceq_{\mu} t^*$. Clearly, if  $s\in {\mathbb B}_1$ then $v_1\in {\mathbb B}_1$, and note that
        if $s\in U$  we have that $I\supsetneqq J(o(s))$ by Claim 3.3(a), and hence $J\supsetneqq [\pi(s^*),\pi(s^*) +2 )$ by $(*)(7)$, and therefore $\pi(s^*) > J^-$ by $(*)(8)$. Analogously, if  $t\in {\mathbb B}_1$ then $v_2\in {\mathbb B}_1$,  and  if $t\in U$  we have that $I\supsetneqq J(o(t))$ by Claim 3.3(b), and hence $K\supsetneqq [\pi(t^*),\pi(t^*) +2 )$ by $(*)(7)$, and therefore $\pi(t^*) > K^-$ by $(*)(8)$. Thus, $v_1,v_2\in {\mathbb B}_1$. Our aim is to show that there is an element $v\in Z$ such that $u\preceq v \prec s,t$ and $\pi(v)\in h\{s,t\}$. We consider the following subcases.

     \vspace{2mm}
      \noindent {\bf Case 2.3.1}. $v_1\in x$.

      \vspace{1mm} Since $u \prec_{\nu} v_1,w$, there is $v'\in i_{\nu}\{v_1,w\}$ such that $u\preceq_{\nu} v'$. So as $v'\in x$ and $v' \prec s^*,t^*$, we can apply  Case 2.2 to $s$, $t$ and $v'$.

      \vspace{2mm}
      \noindent {\bf Case 2.3.2}. $v_2\in x$.

      \vspace{1mm} We see that $s^*$ and $v_2$ are incomparable in $p_{\nu}$, because $s^* \not\prec t^*$. Then as $u \prec_{\nu} s^*, v_2$, there is $v'\in i_{\nu}\{s^*,v_2\}$ such that $u\preceq_{\nu} v'$. Therefore, $v' \prec s^*,t^*$ and $\pi(v')\in h\{s^*,v_2\} = G_I\{J^-,K^-\} \cup \{I^-\}$. Then as $I(\pi(v_2), k+1) = I(\pi(t^*), k+1)$, $v'$ is as required.

      \vspace{2mm}
      \noindent {\bf Case 2.3.3}. $v_1\not\in x$ and $v_2\not\in x$.

      \vspace{1mm} First, suppose that $v_1 \prec_{\nu} w$. Then as
       $v_1 \prec_{\nu} s^*,w$, there is $v'\in i_{\nu}\{s^*,w\}$ with $v_1\preceq_{\nu} v'$. 
       Now as $h\{s^*,w\} = G\{\pi(s^*),\pi(w)\} = G_I\{J^-, \pi(w)\}\cup \{I^-\}$, we infer from $(P3)(d)$ that $\pi(v') < J^-$, and so we have  $\pi(v_1) < J^-$  (because $v_1\preceq_{\nu} v'$), which is impossible.

      \vspace{2mm} Also, note that if  $w \prec_{\nu} v_1$, then $w\prec  s^*,t^*$,  and so we can apply  Case 2.2 to $s$, $t$ and $w$.

      \vspace{2mm}  Now, we assume that $v_1$ and $w$ are incomparable in $p_{\nu}$. Then
       since $u \prec_{\nu} v_1,w$, by $(P3)(f)$ and $(P3)(d)$, we have that there is $v\in i_{\nu}\{v_1,w\}$ such that $u\preceq_{\nu} v$ and
       $\pi(v)\in h\{v_1,w\} = G_I\{J^-,\pi(w)\}\cup \{I^-\}$. So, as $\pi(u) > I^-$ we deduce that $\pi(v)\in G_I\{J^-,\pi(w)\} =  G_I\{\pi(v_1),\pi(w)\}  $.
        As $v \prec_{\nu} v_1,w$, we infer that $v \prec s^*,t^*$. Our aim is to show that $\pi(v)\in h\{s,t\}$. First, note that since $w\prec_{\mu} v_2$ and $\pi(v_2)\in E(I)$,
        we deduce from Claim 3.1 and condition $(C)(8)$ that $\pi(g_{\nu,\nu'}(v_2))\in E(I)$ for every $\nu'\in A$ with $\nu < \nu'$. And since $v_2\in {\mathbb B}_1\setminus x$, we infer that for $\nu'\in A$ with $\nu < \nu'$ we have that $\pi(g_{\nu,\nu'}(v_2))\in {\mathbb B}_1\setminus x$ by conditions $(C1)$ and $(C4)$, and hence $\{\pi(g_{\nu,\nu'}(v_2)):  \nu < \nu', \nu' \in A \}$ is uncountable. Therefore, $\{c_{\gamma}\cap E(I): \gamma \in A\}$ is uncountable.
        We  may assume that $|E(I)| = \om_2$. Otherwise, the argument is simpler. Then as  $\{p_{\gamma} : \gamma \in A \}$ is admissible, $\{c_{\gamma}\cap E(I) : \gamma\in A \}$ is adequate for the function $G_I$. Now since $\pi(w) < \pi(v_2)$, it follows that $G_I\{\pi(v_1),\pi(w)\} \subset G_I\{\pi(v_1),\pi(v_2)\}$, and thus $\pi(v)\in  G_I \{\pi(v_1), \pi(v_2)\}$. Now since $I(\pi(s^*),k+1) = I(\pi(v_1),k+1)$ and $I(\pi(t^*),k+1) = I(\pi(v_2),k+1)$, we infer that $\pi(v)\in h\{s^*,t^*\} = h\{s,t\}$. So, we have finished the proof of  Case 2.3.

       \vspace{2mm} Also, if either $\pi(u)\in a_{\nu}\setminus (a\cup I)$ or $\pi(u)\in a_{\mu}\setminus (a\cup I)$, we can use an argument similar to the one given in Case 2.1. And if $\pi(u)\in (a_{\mu}\setminus a)\cap I$, we can proceed by means of an argument similar to the one given in Case 2.3.

      \vspace{2mm}
       Now, we prove that $p$ satisfies condition $(P4)$. So, consider $s,t\in x_{\nu} \cup x_{\mu}$ such that $s \prec t$ and $s\in {\mathbb B}_1$. If either $s,t\in x_{\nu}$ or  $s,t\in x_{\mu}$, we are done because $p_{\nu}$ and $p_{\mu}$ satisfy $(P4)$. So, assume that $s\in x_{\nu}\setminus x_{\mu}$ and $t\in x_{\mu}\setminus x_{\nu}$. Let $v\in x$ be such that $s \prec_{\nu} v \prec_{\mu} t$. We put $k = \kop(\pi(s),\pi(t))$. Let $J = I(\pi(s),k+1)$ and $K = I(\pi(t),k+1)$. Let $\al = \pi(s)$, $\be = \pi(t)$ and $\gamma =  \pi(v)$. Note that, by condition $(B)(a)-(b)$, we have that $s,v\in {\mathbb B}_1$.  We distinguish the following cases.

        \vspace{2mm}
      \noindent {\bf Case 1}. $\pi(v)\in K$.

      \vspace{1mm} Note that $k = \kop(\pi(s),\pi(v))$. Let $l = \kop(\pi(v),\pi(t))$. We have $l > k$. We assume that $\lop(\al) > k+1$ and $\lop(\gamma) > l+1$. Otherwise, the considerations are simpler. As $s \prec_{\nu} v$, there is a walk $\langle s, s_{\lop(\al)},s_{\lop(\al)-1}, \dots,  s_{k+1}, v_{k+1}, \dots ,$ $v_{\lop(\ga) - 1}, v \rangle$
       from $s$ to $v$ in $p_{\nu}$. And as $v \prec_{\mu} t$, there is a walk $\langle v,v'_{\lop(\ga)}, v'_{\lop(\ga)-1},$ $\dots, v'_{l+1},$ $t_{l+1}, \dots , t \rangle$ from $v$ to $t$ in $p_{\mu}$. Then, the sequence $\langle s, s_{\lop(\al)}, s_{\lop(\al)-1}, \dots,$  $s_{k+1},$ $ v_{k+1}, \dots , v_{l}, t_{l+1}, \dots ,t \rangle$ is a walk from  $s$ to $t$ in $p$.

      \vspace{2mm}
      \noindent {\bf Case 2}. $\pi(v)\in J$.

      \vspace{1mm} Clearly, $k = \kop(\pi(v),\pi(t))$. Let $l = \kop(\pi(s),\pi(v))$. Note that $l > k$.  We  assume that $\lop(\al) > l+1$. Otherwise, the argument is simpler.
      Since $s \prec_{\nu} v$, there is a walk $\langle s, s_{\lop(\al)},$ $s_{\lop(\al)-1}, \dots, s_{l+1}, v'_{l+1}, \dots ,$ $v'_{\lop(\ga) -1}, v \rangle$ from $s$ to $v$ in $p_{\nu}$. And since
      $v \prec_{\mu} t$, there is a walk $\langle v,$ $v_{\lop(\ga)}, v_{\lop(\ga)-1}, \dots, , v_{k+1}, t_{k+1}, \dots , t \rangle$ from $v$ to $t$ in $p_{\mu}$. So, the sequence
      $\langle s, s_{\lop(\al)}, s_{\lop(\al)-1}, \dots, s_{l+1},$ $ v_l, \dots , v_{k+1}, t_{k+1},\dots , t \rangle$ is a walk from $s$ to $t$ in $p$.

      \vspace{2mm}
      \noindent {\bf Case 3}. $\pi(v)\not\in J \cup  K$.

       \vspace{1mm} Clearly, we have $ k = \kop(\pi(s),\pi(v)) = \kop(\pi(v),\pi(t))$. We  suppose that $\lop(\al),\lop(\gamma) > k+1$. 
        Then, if $\langle s, s_{\lop(\al)}, s_{\lop(\al)-1}, \dots,  s_{k+1},$ $ v_{k+1}, \dots , v_{\lop(\ga) - 1}, v \rangle$ is a walk from $s$ to $v$ in $p_{\nu}$ and $\langle v, v'_{\lop(\ga)}, v'_{\lop(\ga)-1}, \dots,  v'_{k+1},$ $t_{k+1}, \dots , t \rangle$ is a walk from $v$ to $t$ in $p_{\mu}$, we have that $\langle s, s_{\lop(\al)}, s_{\lop(\al)-1}, \dots, s_{k+1},$ $t_{k+1}, \dots , t \rangle$ is a walk from $s$ to $t$ in $p$.

       \vspace{2mm} Thus $p\in P_{\delta}$, and so $p\leq_{\de} p_{\nu},p_{\mu}$. This finishes the proof of Lemma 2.10.  $\square$
       
       \newpage {\bf Declaration of competing interest}
       
       \vspace{2mm} The authors declare that they have no known competing financial interests or personal relationships that could have appeared to influence the work reported in this paper.
       
       \vspace{8mm}{\bf Data availability}
       
       \vspace{2mm} No data was used for the research described in the article.


\vspace{12mm}
   \begin{thebibliography}{9}


   \bibitem{ba} J. Bagaria,  Thin-tall spaces and cardinal sequences, in: Open problems in Topology II, Elsevier, Amsterdam, 2007, pp. 115-124.


\bibitem{bs} J. E. Baumgartner, S. Shelah, Remarks on
superatomic Boolean algebras,  Ann. Pure Appl. Log.  33 (2) (l987) 109-129.

\bibitem{bm} M. R. Burke, M. Magidor, Shelah's pcf theory and its applications, Ann. Pure Appl. Log. 50 (3) (1990) 207-254.



\bibitem{ev} K. Er-Rhaimini, B. Veli{\v c}kovi{\' c}, PCF structures of height less than $\om_3$, J. Symb. Log. 75 (4) (2010) 1231-1248.



\bibitem{jsw} I. Juh\'asz, L. Soukup, W. Weiss,  Cardinal sequences of length $< \omega_2$ under GCH, Fundam. Math. 189 (1) (2006)  35-52.

\bibitem{jw} I. Juh\'asz, W. Weiss, Cardinal sequences,  Ann. Pure Appl. Log.  144 (1-3) (2006) 96-106.


\bibitem{ma1} J. C. Mart\'{\i}nez,  A consistency result on thin-very tall Boolean algebras, Isr. J. Math. 123 (2001) 273-284.

\bibitem{ma2} J. C. Mart\'{\i}nez, Cardinal sequences for superatomic Boolean algebras, in: Infinity, computability and metamathematics, in : Tributes, vol 23, College Publications, London, 2014,  pp. 273-284.


\bibitem{ms} J. C. Mart\'{\i}nez, L. Soukup, On cardinal sequences of length $< \om_3$, Topol. Appl. 260 (5) (2019) 116-125.

\bibitem{ru} J. C. Ruyle, Cardinal sequences of PCF structures, Ph. D. thesis, University of California, Irvine, 1998.

\bibitem{se} S. Shelah, Cardinal arithmetic, Oxford Logic Guides, vol. 29, Oxford University Press, 1994.


\bibitem{so} L. Soukup,  Wide scattered spaces and morasses, Topol.  Appl. 158 (5) (2011)  697-707.

\bibitem{to1} S. Todor\v{c}evi\'c,  Partitioning pairs of countable ordinals, Acta Math.  159 (3-4) (1987)  261-294.

\bibitem{to2} S. Todor\v{c}evi\'c,   Walks on ordinals and their characteristics, Progress in Mathematics,
vol. 263,  Birkh\"{a}user Verlag, Basel, 2007.

\end{thebibliography}
\end{document}